

\magnification=\magstep1

\hsize=14cm
\vsize=20cm
\hoffset=-0.3cm
\voffset=0cm
\footline={\hss{\vbox to 1cm{\vfil\hbox{\rm\folio}}}\hss}

\input amssym.def
\input amssym.tex


\font\title=cmr12 at 16pt

\font\teneusm=eusm10
\font\seveneusm=eusm7
\font\fiveeusm=eusm5
\newfam\eusmfam
\def\eusm{\fam\eusmfam\teneusm}
\textfont\eusmfam=\teneusm
\scriptfont\eusmfam=\seveneusm
\scriptscriptfont\eusmfam=\fiveeusm

\font\teneufm=eufm10
\font\seveneufm=eufm7
\font\fiveeufm=eufm5
\newfam\eufmfam

\textfont\eufmfam=\teneufm
\scriptfont\eufmfam=\seveneufm
\scriptscriptfont\eufmfam=\fiveeufm

\def\varGamma{{\mit \Gamma}}
\def\Re{{\rm Re}\,}
\def\Im{{\rm Im}\,}
\def\txt#1{{\textstyle{#1}}}
\def\scr#1{{\scriptstyle{#1}}}
\def\r#1{{\rm #1}}
\def\B#1{{\Bbb #1}}
\def\e#1{{\eusm #1}}
\def\b#1{{\bf #1}}
\def\sgn{{\rm sgn}}

\centerline{\title A new approach to the spectral theory}
\bigskip
\centerline{\title of the fourth moment of the Riemann zeta-function}
\vskip 1cm
\centerline{By {\it Roelof W. Bruggeman} at Utrecht and {\it Yoichi  
Motohashi} at Tokyo}
\vskip 1.5cm {\bf Abstract.}\quad The aim of the present article is
to exhibit a new proof of the explicit formula for the fourth moment
of the Riemann zeta-function that was established by the second
named author a decade   ago. Our proof is new, particularly in that
it dispenses altogether with the spectral theory of sums of
Kloosterman   sums that played a predominant r\^ole in the former
proof. Our argument is, instead, built directly upon the spectral
structure of the space
$L^2(\varGamma\backslash G)$, with
$\varGamma=\r{PSL}_2(\B{Z})$ and $G=\r{PSL}_2(\B{R})$. The
discussion   below thus seems to provide a new insight into the
nature of the Riemann zeta-function, especially in its relation with
automorphic forms over   linear Lie groups that has been perceived
by many.
\par The plan of the paper is as follows: In Section 1 we discuss
salient points of the former proof and describe the explicit formula
in   a conventional fashion. In Section 2 its reformulation is
presented in   terms of automorphic representations occurring in
$L^2(\varGamma\backslash G)$. With this, our motivation is precisely
related. We then proceed to our new proof. In Section 3 we construct
a Poincar\'e series over $G$ whose value at the   unit element is
close to the non-diagonal part of the fourth moment in question. In
Section 4 we develop an account of the Kirillov   scheme, with
which, in Section 5, projections of the Poincar\'e series into  
irreducible subspaces of $L^2(\varGamma\backslash G)$ are explicitly 
calculated in terms of the seed function. Then, in Section 6 a
limiting procedure with respect to the seed is performed, and we
reach a basic spectral expression. This   ends in effect our proof
of the explicit formula, since it remains to appeal   to a process
of analytic continuation, which is, however, the same as the
corresponding part of the former proof, and can largely be omitted.
\par Notations are introduced where they are needed for the first
time, and will continue to be effective thereafter. The parameters
$C>0$ and $\varepsilon>0$ are assumed locally to be
constants arbitrarily large and small, respectively. The
dependency of implicit constants on them can be inferred from the
context.
(Secondary)
\vskip 1.5cm
\centerline {\bf 1. Explicit formula}
\bigskip We formulate the $k$-th moment of the Riemann zeta-function
by
$$
\e{M}_k(\zeta;g)=\int_{-\infty}^\infty
\left|\zeta\left({1\over2}+it\right)\right|^k g(t)dt,
\leqno(1.1)
$$ 
where $k$ is an arbitrary fixed number, and the weight function
$g$ is   assumed, for the sake of simplicity but without much loss
of generality, to be   even, entire, real on $\Bbb R$, and of fast
decay in any fixed horizontal strip. These quantities have been
regarded as major subjects in Analytic Number Theory not only
because they have important applications to a   variety of classical
problems such as the distribution of prime numbers but   also as
they are indispensable means to reveal the intriguing nature of the  
Riemann zeta-function.  Laying stress upon the latter aspect,  we
shall deal with the case $k=4$, the fourth moment, which has
perhaps the richest history of investigations in the theory of mean
values of the zeta-function.
\medskip 
The explicit formula for $\e{M}_4(\zeta;g)$ appeared first
in [17] and later in [19] with certain sophistication. It resulted
from an attempt to generalize the method with which F.V. Atkinson
[1] arrived at an explicit formula for the unweighted mean square. 
Atkinson's formula can indeed be reckoned as the first explicit
result in the theory of mean values of zeta and $L$-functions. In
retrospect, the novelty of his idea rests in that he saw a lattice
structure in a certain non-diagonal expression   which is similar to
$(1.6)$ below. That  probably originated in the works by   H. Weyl
and J.G. van der Corput on the estimation of trigonometrical sums.
To   make his observation effective, Atkinson appealed to two
fundamental implements, analytic continuation and the Poisson sum
formula. Because of its   essential relevance to our present
purpose, we shall summarize Sections 4.1--4.3   of [19] and indicate
how his observation extends to the fourth moment situation:
\medskip Atkinson's argument, when applied to $\e{M}_2(\zeta;g)$,
starts with the introduction of the function
$$
\leqalignno{ I(w_1,w_2;g)&=\int_{-\infty}^\infty
\zeta(w_1-it)\zeta(w_2+it)g(t)dt,&(1.2)\cr
&=\sum_{a=1}^\infty\sum_{b=1}^\infty a^{-w_1}b^{-w_2}
\hat{g}\left({1\over2\pi}\log{a\over b}\right), }
$$ 
where $(w_1,w_2)\in{\Bbb C}^2$ is in the region of absolute
convergence, and
$$
\hat{g}(x)=\int_{-\infty}^\infty g(t)\exp(-2\pi ixt)dt.\leqno(1.3)
$$ 
A shift, either upward or downward, of the contour in
$(1.2)$ gives readily that $I(w_1,w_2;g)$ is meromorphic over the  
entire ${\Bbb C}^2$ and regular in a neighbourhood of the point
$\left({1\over2},{1\over2}\right)$, and that
$$
\e{M}_2(\zeta;g)=I\left({1\over2},{1\over2};g\right) +2\pi\Re
g\left({1\over2}i\right).
\leqno(1.4)
$$ 
Thus $\e{M}_2(\zeta;g)$ can be viewed as a special value of a
meromorphic function of two variables. To attain a continuation of
$I(w_1,w_2;g)$ that could be truly informative, the double  sum in
$(1.2)$ is split into three parts according as $a=b$,  
$a<b$ and $a>b$ so that
$$
I(w_1,w_2;g)=\zeta(w_1+w_2)\hat{g}(0)+I_1(w_1,w_2;g)+I_1(w_2,w_1;g),
\leqno(1.5)
$$ 
where
$$ 
I_1(w_1,w_2;g)=\sum_{n=1}^\infty\sum_{m=1}^\infty{1\over
m^{w_1}(m+n)^{w_2}}\hat{g}
\left({1\over2\pi}\log\left(1+{n\over m}\right)\right).\leqno(1.6)
$$ 
Then $I_1(w_1,w_2;g)$ has to be continued analytically to a  
neighbourhood of
$\left({1\over2},{1\over2}\right)$. This is achieved by regarding
the   inner sum as the one over the lattice $\Bbb Z$. An application
of the Poisson   sum formula yields a meromorphic continuation of
$I_1(w_1,w_2;g)$ to the whole ${\Bbb C}^2$.  Hence, the expression
$(1.5)$ holds throughout  
$\B{C}^2$ as a relation of the four meromorphic functions. The
specialization
$(w_1,w_2)\to\left({1\over2},{1\over2}\right)$ is performed in the  
resulting spectral decomposition of the right side of
$(1.5)$, and  via $(1.4)$ the following explicit formula arises:
$$
\leqalignno{
\e{M}_2(\zeta;g)&=
\int_{-\infty}^\infty\left[\Re{\Gamma'\over\Gamma}\left({1\over2}
+it\right)+2c_E-\log2\pi\right]g(t)dt+2\pi\Re
g\left({1\over2}i\right) &(1.7)\cr &+4\sum_{n=1}^\infty
d(n)\int_0^\infty{1\over\sqrt{r(r+1)}}
g_c\left({1\over2\pi}\log\left(1+{1\over r}\right)\right)
\cos(2\pi{n}r)dr, }
$$ 
where $c_E$ is the Euler constant, $d$ the divisor function, and
$g_c=\Re\hat{g}$ the cosine Fourier transform of $g$ (see [19,  
(4.1.16)]).
\medskip Turning to the fourth moment, we consider the function
$$
\leqalignno{ J(w;g)&=\int_{-\infty}^\infty
\zeta(w_1-it)\zeta(w_2+it)\zeta(w_3+it)
\zeta(w_4-it)g(t)dt&(1.8)\cr
&=\sum_{a=1}^\infty\sum_{b=1}^\infty
\sum_{c=1}^\infty\sum_{d=1}^\infty
a^{-w_1}b^{-w_2}c^{-w_3}d^{-w_4}
\hat{g}\left({1\over2\pi}\log{ad\over bc}\right), }
$$ 
where $w=(w_1,w_2,w_3,w_4)\in{\Bbb C}^4$ is in the region of
absolute convergence. As before, $J(w;g)$ is meromorphic over the
entire ${\Bbb C}^4$, regular in a neighbourhood of the point
$\r{p}_{1\over2}=\left({1\over2},{1\over2},{1\over2},
{1\over2}\right)$,  
and
$$
\e{M}_4(\zeta;g)=J\left(\r{p}_{1\over2};g\right)
-2\pi\Re\left\{\left(c_E-\log2\pi
\right)g\left({1\over2}i\right)+{1\over2}
g'\left({1\over2}i\right)\right \}.
\leqno(1.9)
$$ 
Analogously to $(1.5)$, the quadruple sum in $(1.8)$ is split
into   three parts according as $ad=bc$, $ad<bc$ and $ad>bc$ so that
$$ 
J(w;g)={\zeta(w_1+w_2)\zeta(w_1+w_3)\zeta(w_2+w_4)
\zeta(w_3+w_4)\over\zeta(w_1+w_2+w_3+w_4)}\hat{g}(0)
+J_1(w;g)+J_1(w';g),\leqno(1.10)
$$ 
where $w'=(w_2,w_1,w_4,w_3)$. Our task is then to continue
$J_1(w;g)$ analytically to a neighbourhood of $\r{p}_{1\over2}$.
This is by no means straightforward.  With $J_1(w;g)$, we have
$$ 
J_1(w;g)=\sum_{n=1}^\infty\sum_{m=1}^\infty
{\sigma_{w_1-w_4}(m)\sigma_{w_2-w_3}(m+n)\over
m^{w_1}(m+n)^{w_2}}\hat{g}
\left({1\over2\pi}\log\left(1+{n\over m}\right)\right),\leqno(1.11)
$$ 
where $\sigma_a(n)$ is the sum of $a$-th powers of divisors of
$n$. We face a far more involved expression than $(1.6)$. It is thus
remarkable that Atkinson's view on
$(1.6)$ extends to $(1.11)$. Namely, the inner sum embraces again a
lattice structure. Here the lattice is  
$\varGamma=\r{PSL}_2(\Bbb Z)$, a discrete subgroup of the Lie group
$G=\r{PSL}_2({\Bbb R})$, as to be detailed in the next section.
\par This observation was made in [16] (see also [19, Section 4.2])
but its   proper exploitation was actuated only recently when the
present work was   commenced. In the former proof, instead, the
Ramanujan expansion is employed to   separate the variables $m,n$ in
the arithmetic factor
$\sigma_{w_2-w_3}(m+n)$, and the Vorono{\"\i} scheme is applied to
the arising sum over $m$; namely, the functional equation for the
Estermann zeta-function is invoked. In this way, $J_1(w;g)$ is
transformed into a   sum of  Kloosterman sums. Then the
Kloosterman--Spectral sum formula of N.V. Kuznetsov ([19, Theorems
2.3 and 2.5]) plays a fundamental r\^ole; and a spectral
decomposition of $J_1(w;g)$ emerges. This is the most salient  
point in the former proof, i.e., the assertion in [19, Section 4.5].
The argument is, however, circuitous. In the present paper we aim
to reach the same with   a direct reasoning based on the above
observation.
\par 
With the spectral decomposition of $J_1(w;g)$ thus obtained,
the   argument in [19, Section 4.6] establishes the existence of
$J_1(w;g)$ as a   meromorphic function over the entire $\B{C}^4$.
Hence, the expression $(1.10)$ holds throughout $\B{C}^4$ as a
relation of the four meromorphic functions.   The specialization
$w\to \r{p}_{1\over2}$ gives rise to the explicit   formula for
$\e{M}_4(\zeta;g)$.
\medskip 
We now define the basic spectral terms in the context of
[19]: Let
$\Bbb H$ be the hyperbolic upper half plane
$\left\{z:z=x+iy,x\in{\Bbb   R}, y>0\right\}$ equipped with the
invariant measure
$d\mu(z)=dxdy/y^2$. Let
$\left\{\lambda_j=\kappa_j^2+{1\over4}:
\kappa_j>0, j\ge1\right\}\cup\{0\}$ be the discrete spectrum,
arranged   in non-decreasing order, of the hyperbolic Laplacian
$\Delta=-y^2\left(\partial^2_x+\partial^2_y\right)$ acting on the
Hilbert space $L^2(\varGamma\backslash{\Bbb H},d\mu)$   composed of
all $\varGamma$-automorphic functions that are square integrable
over $\varGamma\backslash{\Bbb H}$ against $d\mu$. We denote by
$\psi_j$ an $L^2$-eigenfunction corresponding to $\lambda_j$, i.e.,
$\Delta\psi_j=\lambda_j\psi_j$. It has the Fourier expansion
$$
\psi_j (z)={\sqrt y}\sum_{\scr{n=-\infty}\atop\scr{n\ne
0}}^\infty\rho_j (n)K_{i\kappa_j}(2\pi |n|y)\exp(2\pi
inx),\leqno(1.12)
$$ 
where $K_\nu$ is the $K$-Bessel function of order $\nu$. We may
assume that the set $\left\{\psi_j:j\ge1\right\}$ forms  an
orthonormal system, and that each $\psi_j$ is a simultaneous
eigenfunction of all Hecke operators. The latter means that we have,
for each positive integer $n$,
$$ 
{1\over {\sqrt
n}}\sum_{\scr{a=1}\atop\scr{ad=n}}^n\sum_{b=1}^{d}\psi_j
\left({{az+b}\over d}\right)= t_j (n)\psi_j (z),\leqno(1.13)
$$ 
with a certain real number $t_j (n)$, and also
$\psi_j (-\bar z)=\epsilon_j \psi_j (z)$, with $\epsilon_j =\pm 1$.
In particular, we have
$\rho_j(1)\not=0$, $\rho_j (n)=\rho_j (1)t_j (n)$,
$\rho_j(-n)=\epsilon_j\rho_j(1) t_j(n)$. Then the Hecke series $H_j
(s)$ associated with $\psi_j$ is defined by
$$ 
H_j(s)=\sum_{n=1}^\infty t_j (n)n^{-s},\leqno(1.14)
$$ 
which converges absolutely for $\Re s>1$, and continues to an
entire   function.
\par 
Each $\psi_j$ is called a real analytic cusp form. We next
introduce holomorphic cusp forms: If $\psi(z)$ is holomorphic
throughout $\Bbb H$,
$\psi(z)(dz)^\ell$ is $\varGamma$-invariant, with a certain fixed
positive integer $\ell$, and $\psi(i\infty)=0$, then  
$\psi$ is called a holomorphic cusp form of weight $2\ell$. The
linear space   consisting of all such functions is a Hilbert space
of finite dimension, in which   the norm of $\psi$ is defined to be
the square root of the integral of  
$|y^\ell
\psi(z)|^2$ over
$\varGamma\backslash{\Bbb H}$ against the measure $d\mu$. Let then
$\left\{\psi_{j,\ell}(z): 1\le j\le\vartheta(\ell)\right\}$ be an  
orthonormal basis of this Hilbert space. We have the Fourier
expansion
$$
\psi_{j,\ell}(z)=\sum_{n=1}^\infty \rho_{j,\ell}(n)n^{\ell-{1\over2}}
\exp(2\pi inz).\leqno(1.15)
$$ 
Similarly to $(1.13)$ we may assume that for any positive integer
$n$
$$ 
{1\over {\sqrt n}}\sum_{\scr{a=1}\atop\scr{ad=n}}^n\left({a\over
d}\right)^\ell\sum_{b=1}^{d}\psi_{j,\ell}
\left({{az+b}\over d}\right)= t_{j,\ell} (n)\psi_{j,\ell}  
(z),\leqno(1.16)
$$ 
with a real number $t_{j,\ell}(n)$ so that $\rho_{j,\ell}(1)\ne0$
and
$\rho_{j,\ell}(n)=\rho_{j,\ell}(1)t_{j,\ell}(n)$. As before, the
Hecke series associated to $\psi_{j,\ell}$ is defined by
$$ 
H_{j,\ell}(s)=\sum_{n=1}^\infty t_{j,\ell}(n)n^{-s},\leqno(1.17)
$$ 
which is again an entire function.
\bigskip With this, we may state the explicit formula for
$\e{M}_4(\zeta;g)$:
\medskip
\noindent {\bf Theorem A.}\quad {\it It holds that
$$
\e{M}_4(\zeta;g)=\big\{\e{M}_{4,0}+{\e{M}}_{4,1}+\e{M}_{4,2}
+\e{M}_{4,3}\big\}(\zeta;g),
\leqno(1.18)
$$ 
where
$$
\leqalignno{
\quad\e{M}_{4,0}(\zeta;g)=\int_{-\infty}^\infty
\sum_{{\scr{p,q,u,v\ge0}}\atop{\scr{pu+qv\le4}}} &
c(p,q,u,v)\,\Re\left\{\left({{\Gamma^{(p)}}
\over\Gamma}\right)^u\left({{\Gamma^{(q)}}\over\Gamma}\right)^v
\left({1\over2}+it\right)\right\}g(t)dt &(1.19)\cr
-2\pi&\Re\left\{(c_E-\log2\pi)g\left({1\over2}i\right)+
{1\over2}ig'\left({1\over2}i\right)\right\},\cr }
$$ 
with effectively computable real absolute constants $c(p,q,u,v);$
and
$$
\leqalignno{
\e{M}_{4,1}(\zeta;g)&=\sum_{j=1}^\infty
\alpha_jH_j\left({1\over2}\right)^3\Lambda(i\kappa_j;g),&(1.20)\cr
\e{M}_{4,2}(\zeta;g)&=\sum_{\ell=1}^\infty
\sum_{j=1}^{\vartheta(\ell)}
\alpha_{j,\ell}H_{j,\ell}
\left({1\over2}\right)^3\Lambda
\left(\ell-{1\over2}\,;g\right),&(1.21)\cr
\e{M}_{4,3}(\zeta;g)&=\int_{(0)}
{{\left(\zeta\left({1\over2}+\nu\right)\zeta\left({1\over2}-\nu\right)
\right)^3}
\over{\zeta(1+2\nu)\zeta(1-2\nu)}}\Lambda(\nu;g) {d\nu\over\pi
i},&(1.22)\cr }
$$ 
with
$$
\alpha_j={|\rho_j(1)|^2\over\cosh\pi\kappa_j},
\qquad \alpha_{j,\ell}={\Gamma(2\ell)\over2^{4\ell-1}\pi^{2\ell+1}}
|\rho_{j,\ell}(1)|^2\,.
\leqno(1.23)
$$ 
The contour in $(1.22)$ is the imaginary axis, and
$$
\leqalignno{ &\Lambda(\nu;g)=\int_0^\infty{1\over\sqrt{r(r+1)}}
g_c\left({1\over2\pi}\log\left(1+{1\over r}\right)
\right)&(1.24)\cr
\times{\Re}&\left\{r^{-{1\over2}-\nu}\left(1-{1\over{\sin\pi  
\nu}}\right) {{\Gamma({1\over2}+\nu)^2}\over{
\Gamma(1+2\nu)}}F\left({1\over2}+\nu,{1\over2}+\nu;1+2\nu;-{1\over r}
\right)
\right\}dr, }
$$ 
with $F$ the hypergeometric function.\/ }
\medskip
\noindent 
This is [19, Theorem 4.2], with a minor change of notation.
The right side of $(1.18)$ has a characteristic pertinent to the
spectral structure of
$L^2(\varGamma\backslash G)$ that is developed in the next section.
It   does not contain any trace of the use of Kloosterman sums.
Then, a problem   comes out: Find a way to reach $(1.18)$ as
directly as possible, especially without recourse to the reduction
to the spectral theory of sums of Kloosterman sums. As has been
indicated above, this is the principal motivation of the present
work. In what follows we shall show an answer to this basic
problem in the theory of the Riemann   zeta-function. It is a
realization of the programme given in [19, Section 4.2].
\bigskip
\noindent {\bf Remark.} A. Ivi\'c's lecture notes [12] give a
thorough account of the theory of mean values of the Riemann
zeta-function. Some major applications of Theorem A are given, and
also in [19]. A recent notable contribution is in Ivi\'c [13]. Those
are, in a   variety of ways, generalizations of consequences derived
from the formula $(1.7)$.   A typical instance is the following
assertion, which could be regarded as an analogue of the localized
version of Atkinson's formula [1]: Let  
$T$ tend to infinity and  
$T^{1\over2}(\log T)^{-C}\le G\le T(\log T)^{-1}$. Then we have
$$
\leqalignno{ {1\over\sqrt{\pi}G}&\int_{-\infty}^\infty
\left|\zeta\left({1\over2}+i(T+t)\right)\right|^4\exp\left(- 
\left({t\over G}\right)^2\right)dt &(1.25)\cr
={\pi\over\sqrt{2T}}\sum_{j=1}^\infty\alpha_j&H_j
\left({1\over2}\right)^3
\kappa_j^{-{1\over2}}\sin\left(\kappa_j\log{\kappa_j\over4eT}\right)
\exp\left(-{1\over4}\left({G\kappa_j\over T}\right)^2\right)
+O\left((\log T)^{3C+9}\right) }
$$ 
(see [19, (5.1.44)]). This makes it clear that the values of the
zeta-function on the   critical line are related to eigenvalues of
the hyperbolic Laplacian.
\vskip 1cm
\centerline{\bf 2. Reformulation}
\bigskip 
Now, we make precise the lattice structure of the function
$J_1(w;g)$: We define a function $g_*$ on $G$ by
$$ g_*(\r{g})= |a|^{-w_1}|b|^{-w_2}|c|^{-w_3}|d|^{-w_4}
\hat{g}\left({1\over2\pi}\log{ad\over bc}\right)\iota(ad),\quad
G\ni\r{g}=\left[\matrix{a&b\cr c& d}\right],\leqno(2.1)
$$ 
where $\iota$ is the characteristic function of the negative
reals, and   the matrix is in the projective sense. Then a
rearrangement gives
$$ J_1(w;g)={1\over4}\sum_{\scr{\r{g}\in
\r{M}_2({\Bbb Z})}\atop\scr{\det\r{g}>0}}
(\det\r{g})^{-z_1-{1\over2}}
g_*\left({\r{g}\over\sqrt{\det\r{g}}}\right),\leqno(2.2)
$$ 
where $w$ is in the region of absolute convergence, and
$$ 
z_1={1\over2}(w_1+w_2+w_3+w_4-1).\leqno(2.3)
$$ 
Invoking Hecke's representatives for the quotient $\r{SL}_2({\Bbb
Z})\backslash\r{M}_2({\Bbb Z})$, we have, in place of $(2.2)$,
$$ 
J_1(w;g)={1\over2}\sum_{n=1}^\infty n^{-z_1-{1\over2}}
\sum_{d|n}\sum_{b\bmod d}\,\sum_{\gamma\in\varGamma}
g_*\left(\gamma\left[\matrix{1&b/d\cr&1}\right]
\left[\matrix{\sqrt{n}/d& \cr&d/\sqrt{n}}\right]\right).\leqno(2.4)
$$ 
The lattice structure is now evident.
\par 
Further, if we put
$$
\e{P}g_*(\r{g})=\sum_{\gamma\in\varGamma} g_*(\gamma\r{g}),\quad
\r{g}\in G,
\leqno(2.5)
$$ 
and
$$
\e{T}=\sum_{n=1}^\infty T_n n^{-z_1}\,,\leqno(2.6)
$$ 
with the Hecke operators $T_n$ which act from the left (cf.
$(1.13)$ and
$(1.16)$). Then $(2.4)$ gives formally
$$ 
J_1(w;g)={1\over2}\e{T}\e{P}g_*(1),\leqno(2.7)
$$ 
where the argument on the right side is the unit element of $G$.
\par 
This is revealing. However, the Poincar\'e series $\e{P}g_*$ is
not   really defined throughout $G$. In fact, the sum $(2.5)$
diverges on a dense   subset of $G$; see the discussion leading
$(3.5)$ below. Nevertheless, $(2.7)$   itself is correct and
demonstrates that
$J_1(w;g)$ is an object closely related to the space of
$\varGamma$-automorphic functions over $G$.
\medskip We now start pursuing this line of reasoning.
\medskip To begin with, we collect here elements of the theory of
$\varGamma$-automorphic representations of $G$:  We write
$$
\r{n}[x]=\left[\matrix{1&x\cr&1}\right],\quad
\r{a}[y]=\left[\matrix{\sqrt{y}&\cr&1/\sqrt{y}}\right],\quad
\r{k}[\theta]=\left[\matrix{\phantom{-}\cos\theta&
\sin\theta\cr-\sin\theta &\cos\theta}\right].\leqno(2.8)
$$ 
Let $N=\left\{\r{n}[x]: x\in\B{R}\right\}$, $A=\left\{\r{a}[y]:  
y>0\right\}$, and
$K=\left\{\r{k}[\theta]:\theta\in\B{R}/\pi\B{Z}\right\}$ so that
$G=NAK$ be the Iwasawa decomposition of the Lie group $G$.  We  
read it as $G\ni
\r{g}=\r{n}\r{a}\r{k}=\r{n}[x]\r{a}[y]\r{k}[\theta]$;   throughout
the sequel, the coordinate $(x,y,\theta)$ retains this definition.
The Haar measures on the groups
$N$, $A$, $K$, $G$ are defined, respectively, by $d\r{n}=dx$,
$d\r{a}=dy/y$, $d\r{k}=d\theta/\pi$,
$d\r{g}=d\r{n}d\r{a}d\r{k}/y$,  with Lebesgue measures $dx$,
$dy$, $d\theta$. The Lie algebra of $G$ is spanned by
$$ 
{\bf X}=\pmatrix{0&1\cr0&0},\quad {\bf Y}=\pmatrix{0&0\cr1&0},
\quad{\bf H}=\pmatrix{1&\phantom{-}0\cr0&-1}.\leqno(2.9)
$$ 
The universal enveloping algebra of $G$ is denoted by $\e{U}$.
Its   center is the polynomial ring on the Casimir element  
$\Omega=y^2\left(\partial_x^2+
\partial_y^2\right)-y\partial_x
\partial_\theta$.
\medskip The space $L^2(\varGamma\backslash G)$ is composed of all
left
$\varGamma$-automorphic functions on $G$, vectors for short, which
are square integrable over $\varGamma\backslash G$ against $d\r{g}$.
Elements of $G$ act unitarily on vectors from the right. We have the
orthogonal decomposition into invariant subspaces
$$ 
L^2(\varGamma\backslash G)=\B{C}\cdot1\bigoplus
{}^0\!L^2(\varGamma\backslash G)\bigoplus
{}^e\!L^2(\varGamma\backslash G).\leqno(2.10)
$$ 
Here ${}^0\!L^2$ is the cuspidal subspace spanned by vectors
whose Fourier expansions with respect to the left action of $N$ have
vanishing constant terms. The subspace
${}^e\!L^2$ is spanned by integrals of Eisenstein series, as is to be
detailed in Lemma 2 below. Note that invariant subspaces and
$\varGamma$-automorphic representations of $G$ are interchangeable  
concepts, and we refer to them in a mixed way.
\par The cuspidal subspace splits into irreducible subspaces:
$$ 
{}^0\!L^2(\varGamma\backslash G)=\overline{\bigoplus V}.
\leqno(2.11)
$$ 
The Casimir operator becomes a constant multiplication in each
$V$; that is,
$$
\Omega|_{V^\infty}=\left(\nu_V^2-{1\over4}\right)\cdot1,\leqno(2.12)
$$ 
where $V^\infty$ is the set of all smooth vectors in $V$. Under
our present supposition that
$\varGamma=\r{PSL}_2(\B{Z})$, we can restrict our attention to two
cases: either $i\nu_V<0$ or $\nu_V$ is equal to half a positive odd  
integer. According to the right action of $K$, the space
$V$ is decomposed into $K$-irreducible subspaces
$$ 
V=\overline{\bigoplus_{p=-\infty}^\infty V_p},\quad \dim V_p\le
1.\leqno(2.13)
$$ 
If it is not trivial, $V_p$ is spanned by a
$\varGamma$-automorphic function on which the right translation by
$\r{k}[\theta]$ becomes the multiplication by the factor
$\exp(2ip\theta)$. It is called a $\varGamma$-automorphic form of
spectral parameter $\nu_V$ and weight $2p$.
\medskip Let us assume temporarily that $V$ belongs to the unitary
principal series, i.e., $i\nu_V<0$, under our present situation. 
Then one can   show that
$\dim V_p=1$ for all
$p\in\B{Z}$ and that there exists a complete orthonormal system
$\{\varphi_p\in V_p:\,p\in\B{Z}\}$ of $V$ such that
$$
\varphi_p(\r{g})=\sum_{\scr{n=-\infty}\atop\scr{n\ne0}}
^\infty{\varrho_V(n)\over\sqrt{|n|}}
\e{A}^{\sgn(n)}\phi_p(\r{a}[|n|]\r{g};\nu_V),
\leqno(2.14)
$$ 
where $\phi_p(\r{g};\nu)=y^{\nu+{1\over2}}\exp(2ip\theta)$, and
$$
\e{A}^\delta\phi_p(\r{g};\nu)=\int_{-\infty}^\infty
\exp(-2\pi i\delta x)\phi_p (\r{w}\r{n}[x]\r{g};\nu)dx,\quad
\r{w}=\left[\matrix{&1\cr-1&\cr}
\right].\leqno(2.15)
$$ 
The $\e{A}^\delta$ is a specialization of the Jacquet operator.
This follows from a study of the Fourier expansion of $\varphi_p$  
coupled with the action of the Maass operators
${\bf E}= e^{2i\theta}(2iy\partial_x+2y\partial_y-i\partial_\theta)$
and
$\overline{\bf E}$. It should be observed that the coefficients
$\varrho_V(n)$ in $(2.14)$ do not depend on the weight.
\par We note that
$$
\leqalignno{
\e{A}^\delta\phi_p(\r{g};\nu)&=y^{{1\over2}-\nu}\exp(2\pi i\delta x)
\int_{-\infty}^\infty{\exp(2\pi iy\xi)
\over(\xi^2+1)^{{1\over2}+\nu}}
\left({\xi+i\over\xi-i}\right)^{\delta  
p}d\xi\cdot\exp(2pi\theta)&(2.16)\cr &=(-1)^p\pi^{{1\over2}+\nu}
\exp(2\pi i\delta x) {W_{\delta p,\nu} (4\pi y)\over\Gamma(\delta
p+{1\over2}+\nu)}\exp(2pi\theta), }
$$ 
where $W_{\lambda,\mu}(y)$ is the Whittaker function (see [24,
Chapter   XVI]). The first line is valid for
$\Re\nu>0$, while the second defines $\e{A}^\delta\phi_p$ for all
$\nu\in\B{C}$. In particular, we have the expansion
$$
\varphi_0(\r{g})={2\pi^{{1\over2}+\nu_V}\over\Gamma({1\over2}+\nu_V)}
\sqrt{y}\sum_{\scr{n=-\infty}\atop\scr{n\ne0}}
^\infty\varrho_V(n)K_{\nu_V}(2\pi|n|y)
\exp(2\pi inx).\leqno(2.17)
$$ 
This corresponds to $(1.12)$. Namely, on the identification
$$
\varphi_0(\r{g})=\psi_j(x+iy),\leqno(2.18)
$$ 
we have
$$
\nu_V=i\kappa_j,\quad
\varrho_V(n)={\Gamma({1\over2}+\nu_V)\over2\pi^{{1\over2}+\nu_V}}
\rho_j(n).\leqno(2.19)
$$
\medskip 
Next, let us consider a $V$ in the discrete series; that
is, $\nu_V=\ell-{1\over2}$, $1\le\ell\in\B{Z}$. We have, in place  
of
$(2.13)$,
$$
\hbox{either\quad
$\displaystyle{V=\overline{\bigoplus_{p=\ell}^\infty   V_p}}$
\quad or\quad
$\displaystyle{V=\overline{\bigoplus_{p=-\infty}^{-\ell} V_p}}$}\,,
\leqno(2.20)
$$ 
with $\dim V_p=1$, corresponding to the holomorphic and the
antiholomorphic discrete series. The involution
$\omega:\,\r{g}=\r{nak}\mapsto\r{n}^{-1}\r{a}\r{k}^{-1}$ maps one to
the other.  In the holomorphic case, we have a complete orthonormal
system $\{\varphi_p:\, p\ge\ell\}$ in $V$ such that
$$
\varphi_p(\r{g})=\pi^{{1\over2}-\ell}\left({\Gamma(p+\ell)\over
\Gamma(p-\ell+1)}\right)^{1\over2}\sum_{n=1}^\infty
{\varrho_V(n)\over\sqrt{n}}
\e{A}^+\phi_p\left(\r{a}[n]\r{g};\ell-{1\over2}\right),\leqno(2.21)
$$ 
In particular, we have
$$
\varphi_\ell(\r{g})=(-1)^\ell{2^{2\ell}\pi^{\ell+{1\over2}}\over
\sqrt{\Gamma(2\ell)}}\exp(2i\ell\theta)y^\ell
\sum_{n=1}^\infty\varrho_V(n)n^{\ell-{1\over2}}\exp(2\pi in(x+iy)),
\leqno(2.22)
$$ 
which corresponds to $(1.15)$. Thus, on the identification
$$
\varphi_\ell(\r{g})=y^\ell\exp(2i\ell\theta)\psi_{j,\ell}(x+iy),
\leqno(2.23)
$$ 
we have
$$
\nu_V=\ell-{1\over2},\quad  
\varrho_V(n)=(-1)^\ell{\sqrt{\Gamma(2\ell)}\over
2^{2\ell}\pi^{\ell+{1\over2}}}\rho_{j,\ell}(n).\leqno(2.24)
$$ 
On the other hand, if $V$ is in the antiholomorphic discrete
series,   then we have a complete orthonormal system $\{\varphi_p:\,
p\le-\ell\}$ in $V$ such that
$$
\varphi_p(\r{g})=\pi^{{1\over2}-\ell}\left({\Gamma(|p|+\ell)\over
\Gamma(|p|-\ell+1)}\right)^{1\over2}\sum_{n=-\infty}^{-1}
{\varrho_V(n)\over\sqrt{|n|}}
\e{A}^-\phi_p\left(\r{a}[|n|]\r{g};\ell-{1\over2}\right),\leqno(2.25)
$$ 
with $\varrho_V(n)=\varrho_{V^*}(-n)$, where $V^*=\left\{
\varphi\omega:\varphi\in V\right\}$ is in the holomorphic discrete
series.
\medskip 
With this, $(1.13)$ and $(1.16)$ are extended to
$$ 
T_n|_V=t_V(n)\cdot1,\leqno(2.26)
$$ 
for any $V$, where $t_V(n)$ is newly defined to be equal to
either $t_j(n)$ or $t_{j,\ell}(n)$, as specified above. Thus,   for
any non-zero integer $n$,
$$
\varrho_V(n)=\varrho_V(\sgn(n))t_V(|n|).\leqno(2.27)
$$ 
We have here the obvious convention that $\varrho_V(-1)=0$ and
$\varrho_V(1)=0$ for $V$ in the holomorphic and antiholomorphic
discrete series, respectively. If $V$ is in the unitary principal
series, then
$\varrho_V(-1)=\epsilon_V\varrho_V(1)$ following the assertion
adjacent
$(1.13)$, but with $\epsilon_j$ being denoted by
$\epsilon_V$.
\medskip We observe that $(2.19)$, $(2.24)$, and $(2.27)$ translate
$(1.23)$ into
$$
\hbox{$
|\varrho_V(1)|^2+|\varrho_V(-1)|^2
=\displaystyle{1\over2}\alpha_j$\quad
or\quad$\displaystyle{1\over2}\alpha_{j,\ell}$}\,,\leqno(2.28)
$$ 
according to the series to which $V$ belongs. Lemmas 2.3 and 2.4
of [19] yield
$$
\sum_{\scr{V}\atop\scr{|\nu_V|\le N}} |\varrho_V(\pm1)|^2\ll
N^2,\leqno(2.29)
$$ 
as $N$ tends to infinity, with the implied constant being
absolute. Those lemmas give also
$$ 
t_V(n)\ll n^{{1\over4}+\varepsilon}.\leqno(2.30)
$$ 
See the remark at the end of the present section.
\medskip Replacing $(1.14)$ and $(1.17)$, we associate the Hecke
series
$$ 
H_V(s)=\sum_{n=1}^\infty t_V(n)n^{-s}\leqno(2.31)
$$ 
to each $V$, in both types of cuspidal representations. This
continues to an entire function, satisfying the functional equation
$$
\leqalignno{
H_V(s)=&2^{2s-1}\pi^{2(s-1)}\Gamma(1-s+\nu_V)
\Gamma(1-s-\nu_V)&(2.32)\cr
&\times\left\{\epsilon_V\cos\pi\nu_V-\cos\pi s\right\}H_V(1-s), }
$$ 
where $\epsilon_V$ is left undefined for $V$ in the discrete
series,   since
$\cos\pi\nu_V=0$ there. In particular, $H_V(s)$ is of polynomial
growth   in both $s$ and $\nu_V$ in any fixed vertical strip in the
$s$-plane (see   [19, Chapter 3]).
\medskip 
Now, let us reformulate the explicit formula $(1.18)$ in
terms of
$\varGamma$-automorphic representations of $G$. To this end, we  
introduce the Bessel function of representations of $G$, in the
sense of [10]:
$$ 
j_\nu(u)=\pi{\sqrt{|u|}\over\sin\pi\nu}\left(J_{-2\nu}^{\sgn(u)}
(4\pi\sqrt{|u|})-J^{\sgn(u)}_{2\nu}(4\pi\sqrt{|u|})\right),
\leqno(2.33)
$$ where $J^+_\nu=J_\nu$ and $J^-_\nu=I_\nu$ with the ordinary
notation for Bessel functions. Also we put
$$
\Theta(\nu;g)=\int_0^\infty
{1\over\sqrt{r(r+1)}}g_c\left({1\over2\pi}
\log\left(1+{1\over r}\right)\right)
\Xi(r;\nu)dr,\leqno(2.34)
$$ 
with
$$
\Xi(r;\nu)=\int_{\B{R}^\times} j_0(-u)j_\nu\left({u\over r}\right)
{d^\times u\over\sqrt{|u|}},\quad d^\times u={du\over
|u|}.\leqno(2.35)
$$
\medskip 
Then Theorem A takes the following new form:
\medskip
\noindent 
{\bf Theorem B.}\quad {\it We have
$$
\e{M}_4(\zeta;g)=\left\{\e{M}_4^{(r)}+\e{M}_4^{(c)}
+\e{M}_4^{(e)}\right\}(\zeta;g).\leqno(2.36)
$$ 
Here $\e{M}^{(r)}_4=\e{M}_{4,0}$, and
$$
\leqalignno{
\e{M}_4^{(c)}(\zeta;g)&=\sum_{V}
\left(|\varrho_V(1)|^2+|\varrho_V(- 
1)|^2\right)H_V\left({1\over2}\right)^3
\Theta(\nu_V;g),&(2.37)\cr
\e{M}_4^{(e)}(\zeta;g)&=\int_{(0)}
{\left(\zeta\left({1\over2}+\nu\right)
\zeta\left({1\over2}-\nu\right)\right)^3
\over\zeta(1+2\nu)\zeta(1-2\nu)}\Theta(\nu;g){d\nu\over2\pi
i}\,.&(2.38) }
$$ 
The variable $V$ runs over a maximal orthogonal system of
Hecke-invariant cuspidal $\varGamma$-automorphic representations of
$G$. }
\medskip
\noindent 
{\it Proof.\/} This is a corrected version of the
reformulation made in [20] without details. We shall afterward prove
this directly using the spectral expansion of smooth vectors, i.e.,
Lemma 2 below, and the Kirillov scheme developed in Section 4. The
latter will reveal the reason why the integral transform $(2.34)$
comes up in the explicit   formula. Here we show briefly how
$(2.36)$ follows from $(1.18)$. This might   appear redundant but
seems to have its own interest. Also it provides us with   an
opportunity to make certain preparation for our later discussion.
\medskip
We observe two Mellin transforms: If $|\Re\nu|-{1\over2}<\Re
s$, then
$$
\int_{-\infty}^0 j_\nu(u)|u|^{s-1}du={1\over\pi}(2\pi)^{-2s}\cos\pi
\nu\,
\Gamma\left(s+{1\over2}+\nu\right)
\Gamma\left(s+{1\over2}-\nu\right);\leqno(2.39)
$$ 
and if $|\Re\nu|-{1\over2}<\Re s<-{1\over4}$, then
$$
\int_0^\infty j_\nu(u)u^{s-1}du=-{1\over\pi}(2\pi)^{-2s}\sin\pi s\,
\Gamma\left(s+{1\over2}+\nu\right)
\Gamma\left(s+{1\over2}-\nu\right);\leqno(2.40)
$$ 
The former follows from $(8)$ of [23, Section 13.21], and the
latter from $(1)$ of [23, Section 13.24]. Both integrals are
absolutely   convergent in the respective ranges, because of
$(2.45)$ below.
\par To compute the following integral, we replace the factor
$j_\nu(u/r)$ using the Mellin inversion of $(2.39)$ with a  contour
$(\beta)$, i.e., the vertical line $\Re s=\beta$. It is to be chosen
so that the resulting double integral is absolutely convergent. Then
we may exchange the order of integration and use
$(2.40)$ for the inner integral which is a Mellin transform of
$j_0(u)$,
$u>0$. In this way, we have, for
$r>0$ and, e.g., for $-{1\over 4}<\beta<-|\Re\nu|$,
$$
\leqalignno{
\int_{-\infty}^0 &j_0(-u)j_\nu\left({u\over r}\right) {du\over
|u|^{{3\over2}}}&(2.41)\cr
=&{\cos(\pi\nu)\over\pi^2i}\int_{({\beta})}
\cos\pi s\,\Gamma\left(s+{1\over2}+\nu\right)
\Gamma\left(s+{1\over2}-\nu\right)
\Gamma\left(-s\right)^2r^s ds\cr
=&r^{-{1\over2}+\nu}{\Gamma({1\over2}-\nu)^2\over
\Gamma(1-2\nu)}F\left({1\over2}-\nu,{1\over2}-\nu; 1-2\nu;
-{1\over r}\right)
\cr &+\hbox{the same expression but with $\nu\mapsto -\nu$}. 
}
$$ 
For the second equality see e.g. [19, pp.\ 119--120]. Similarly
but exchanging the r\^oles of the factors $j_0(-u)$ and
$j_\nu (u/r)$, we have, on the same condition,
$$
\leqalignno{
\int_0^\infty &j_0(-u)j_\nu\left({u\over r}\right) {du\over
u^{{3\over2}}}&(2.42)\cr =&{1\over\pi^2i}\int_{(\beta)}\cos\pi
s\,\Gamma(-s+\nu)
\Gamma(-s-\nu)\Gamma\left(s+{1\over2}\right)^2r^{-s-{1\over2}}ds\cr
=& {r^{-{1\over2}+\nu}\over\sin\pi\nu}{\Gamma({1\over2}-\nu)^2\over
\Gamma(1-2\nu)}F\left({1\over2}-\nu,{1\over2}-\nu;1-2\nu;-{1\over
r}\right)\cr &+\hbox{the same expression but with $\nu\mapsto -\nu$}.
}
$$
\par 
Thus we find that $\Theta(\nu;g)=2\Lambda(\nu;g)$ for $\nu\in
i{\Bbb R}$. Then $(2.38)$ is immediate. Also, taking into account
the first   case of
$(2.28)$ we obtain $(2.37)$ if $V$ is restricted to the unitary  
principal series.
\par On the other hand, if
$\ell$ is a positive integer, then we have
$j_{\ell-{1\over2}}(u)=0$   for $u<0$ and  
$j_{\ell-{1\over2}}(u)=(-1)^\ell2\pi\sqrt{u}J_{2\ell-1}
(4\pi\sqrt{u})$
for
$u>0$. This gives, for $-\ell<\Re{s}<-{1\over4}$,
$$
\int_0^\infty j_{\ell-{1\over2}}(u)u^{s-1}du=(-1)^\ell
(2\pi)^{-2s}{\Gamma(s+\ell)
\over\Gamma(\ell-s)},\leqno(2.43)
$$ 
and
$$
\int_0^\infty j_0(-u)j_{\ell-{1\over2}}\left({u\over r}\right)
{du\over u^{{3\over2}}}=2 (-r)^{-\ell}{\Gamma(\ell)^2\over
\Gamma(2\ell)}F\left(\ell,\ell;2\ell;-{1\over r}\right).\leqno(2.44)
$$ 
Thus we have
$\Theta\left(\ell-{1\over2};g\right)
=\Lambda\left(\ell-{1\over2};g\right)$ if $\ell$ is even. However,
the   same does not hold for all $\ell$, since  
$\Lambda\left(\ell-{1\over2};g\right)=0$ if $\ell$ is odd.
Nevertheless, $(2.32)$ gives $H_V({1\over2})=0$ if
$\nu_V=2\ell+{1\over2}$. That is, such $V$ are irrelevant. With
this, on noting
$(2.20)$ and the second case of $(2.28)$, we end the proof.
\medskip It should be remarked that if $\nu$ is bounded and
$|\Re\nu|<{1\over2}$, then
$$ 
j_\nu(u)\ll
\cases{
\hfill u^{{1\over4}}\hfill & if $u>1$,\cr
\hfill |u|^{-{1\over \varepsilon}}\hfill & if $u<-1,$\cr
|u|^{{1\over2}-|\Re\nu|-\varepsilon}& if
$|u|\le1$. }\leqno(2.45)
$$ 
Also, for integral $\ell\ge1$
$$ 
j_{\ell-{1\over2}}(u)\ll\min(u^{1\over4},u^\ell),
\quad u>0.\leqno(2.46)
$$ 
In fact, if $|u|$ is small, both follow from power series
expansions of Bessel functions. If $u$ is   negative, then
$j_\nu(u)=4\sqrt{|u|}\cos\pi\nu\, K_{2\nu}(4\pi\sqrt{|u|})$, which
decays exponentially. If $u$ tends to
$+\infty$, $(2.46)$ is a consequence of a well-known asymptotic  
expansion for
$J$-Bessel functions (see [23, Section 7.21]). To treat the
remaining case, we note that the formula
$(12)$ on [23, p.\ 180] gives that for $u>0$,
$|\Re\nu|<{1\over2}$,
$$ 
j_\nu(u)=4\sqrt{u}\int_0^\infty  
\cos(4\pi\sqrt{u}\cosh\xi)\cosh(2\nu\xi)d\xi.
\leqno(2.47)
$$ 
If $u$ is large, then we divide the integral at
$\xi=u^{-{1\over4}}$. We bound the integrand trivially for smaller  
$\xi$; otherwise we integrate in part with respect to the cosine
factor. This   ends the proof of
$(2.45)$--$(2.46)$.
\bigskip
\noindent 
{\bf Remark.} The proofs of $(2.29)$--$(2.30)$ in [19]
might appear to   come rather close to the spectral theory of sums
of Kloosterman sums. A   closer examination will, however, reveal
that the proofs depend only on a non-trivial bound for individual
Kloosterman sums. T. Estermann's elementary bound should work fine
for $(2.29)$. As to $(2.30)$, A.   Weil's bound is used, but it can
be replaced by Estermann's, too, although a   weaker bound results.
In fact, the exponent in $(2.30)$ is not much relevant   to our
discussion. It should be worth stressing that our proof, developed  
below, of Theorem B depends on the theory of Kloosterman sums solely
via $(2.29)$--$(2.30)$.
\par 
The term $\e{M}_{4,2}$ in $(1.18)$ is there because the
Kloosterman-Spectral sum formula is employed in the former proof and
it contains a contribution coming from holomorphic cusp forms. In
$(2.36)$ this part is replaced by the contribution of irreducible
representations belonging to the discrete series. Concerning the sum
formula, the notes [3] by the first named author arose from the wish
to understand various terms by means of the theory of automorphic
representations; see also [4]. In this respect, our present work
shares motivations with [3] to a considerable extent. Besides, we
have been inspired by J.W. Cogdell and I. Piatetski-Shapiro [8].
\vskip 1cm
\centerline{\bf 3. Poincar\'e series}
\bigskip 
We now start a proof of Theorem B by appealing to the
spectral theory of $L^2(\varGamma\backslash G)$. In this section we
fix the Poincar\'e series on which our discussion is to be developed.
\medskip In place of $(2.1)$ we put
$$ f_{\psi\tau}(\r{g})= |a|^{-w_1}|b|^{-w_2}|c|^{-w_3}|d|^{-w_4}
\psi\left({ad\over bc}\right)\tau(ad),\quad
G\ni\r{g}=\left[\matrix{a&b\cr c&d}\right],\leqno(3.1)
$$ 
where
$$
\psi^{(l)}(x)\ll\min\left(|x|^B,|x|^{-B}\right),\leqno(3.2)
$$ 
and
$$
\hbox{$\tau(x)=0,\;x>0$;\quad $\tau^{(l)}(x)\ll\min\left(|x|^B,
|x|^{-B}\right)$,}\leqno(3.3)
$$ for each $l$ and any constant $B>0$. In other words, it is
assumed that   all derivatives of $\psi(x)$ and $\tau(x)$ decay
faster than any power of  
$|x|$ as $x$ tends to $0$ and infinity. The specialization
$\psi(x)=\hat{g}\left({1\over2\pi}\log |x|\right)$ naturally
satisfies $(3.2)$. We consider the 
Poincar\'e series $\e{P}f_{\psi\tau}(\r{g})$ following
$(2.5)$. Ignoring the convergence issue temporarily, we have, as  
$(2.7)$,
$$
{1\over2}\e{T}\e{P}f_{\psi\tau}(1)=
\sum_{n=1}^\infty\sum_{m=1}^\infty
{\sigma_{w_1-w_4}(m)\sigma_{w_2-w_3}(m+n)\over
m^{w_1}(m+n)^{w_2}}\psi\left({m\over m+n}\right)
\tau\left(-{m\over n}\right).\leqno(3.4)
$$ 
One may take the limit as $\tau$ tends to $\iota$ the
characteristic   function of the negative reals. Then the result is
indeed comparable with  
$(2.7)$ (see also $(1.11)$). However, in general the series
$\e{P}f_{\psi\tau}(\r{g})$ does not converge for all $\r{g}$. To
see   this, let $\gamma_0\in\varGamma$ be a hyperbolic element, and
$\r{g}_0\in G$   be such that
$\r{g}_0^{-1}\gamma_0\r{g}_0 =\r{a}[\lambda]$ with $\lambda>1$.
Then,   for any integer $n$, we have
$f_{\psi\tau}(\gamma_0^n\r{g}_0)=f_{\psi\tau}(\r{g}_0)
\lambda^{{1\over2}n(w_2+w_4-w_1-w_3)}$, which obviously implies the
divergence of $\e{P}f_{\psi\tau}$ at
$\r{g}_0$, provided $f_{\psi\tau}(\r{g}_0)\ne0$. Hence
$\e{T}\e{P}f_{\psi\tau}$ is not well-defined. To overcome this
difficulty, we introduce the modification
$$ 
f_{\psi\tau\eta}=|a|^{-w_1}|b|^{-w_2}|c|^{-w_3}|d|^{-w_4}
\psi\left({ad\over bc}\right)\tau(ad)
\eta\left({d\over c}\right),\leqno(3.5)
$$ 
with an $\eta$ satisfying
$$
\eta(-x)=\eta(x);\quad
\eta^{(l)}(x)\ll \min\left(|x|^B,|x|^{-B}\right),\leqno(3.6)
$$ 
as before. Note that
$f_{\psi\tau\eta}$ is left
$A$-equivariant:
$$ 
f_{\psi\tau\eta}(\r{a}[y]\r{g})
=y^{z_2-{1\over2}}f_{\psi\tau\eta}(\r{g}),\quad
z_2={1\over2}(w_3+w_4-w_1-w_2+1).\leqno(3.7)
$$
\medskip
\noindent {\bf Lemma 1.} {\it Let $f=f_{\psi\tau\eta}$, with bounded
$w$ in the domain
$$
\Re(w_3+w_4)>2+\Re(w_1+w_2)>4.\leqno(3.8)
$$ 
Then the Poincar\'e series
$\e{P}f$ converges absolutely and uniformly to an infinitely  
differentiable
$\varGamma$-automorphic function on
$G$.  The same holds for $\e{T}\e{P}f$, and in particular
$$
{1\over2}\e{T}\e{P}f(1)=\sum_{n=1}^\infty\sum_{m=1}^\infty
\psi\left({m\over m+n}\right)\tau\left(-{m\over  
n}\right)\sum_{\scr{ad=m}\atop\scr{bc=m+n}}
a^{-w_1}b^{-w_2}c^{-w_3}d^{-w_4}\eta\left({d\over
c}\right),\leqno(3.9)
$$ 
with positive integers $a,b,c,d$. }
\medskip
\noindent 
{\bf Remark.} Throughout the sequel, this definition for
$f$ will be retained. The condition
$\Re(w_3+w_4)>2+\Re(w_1+w_2)>3$ is sufficient for the convergence of
$\e{P}f$, but $\e{T}\e{P}f$ requires $(3.8)$. The heuristic identity
$(3.4)$ can be understood to be the limit of
$(3.9)$ as $\eta$ tends to the characteristic function of
$\B{R}^\times$. Also, the limiting procedure mentioned above with
respect to $\tau$ is to be considered. We shall perform these with
an explicit choice of $\tau$ and $\eta$, in the final section. The
infinite differentiability is required because of our   appeal to
Lemma 2 below. However, one may take a conceptually simpler
approach,   and then this stringent condition can be avoided, at the
cost of an extra estimation procedure. See the remark at the end of
this section. In particular, the smoothness condition on
$\psi$, $\tau$, $\eta$ introduced above could considerably be
relaxed, although it is irrelevant to our present purpose.
\medskip
\noindent 
{\it Proof.\/} Let $G=AN\sqcup AN\r{w}N$ be the Bruhat
decomposition of $G$, where $\r{w}$, $A$, $N$ are defined above. We
have, for $\sin\theta\not=0$, i.e., in the big cell,
$$
\r{n}[x]\r{a}[y]\r{k}[\theta]=
\left[\matrix{\displaystyle{\sqrt{y}\over\sin\theta}&\cr
&\displaystyle{\sin\theta\over\sqrt{y}}}\right]
\r{n}\left[{x\over y}\sin^2\theta
-\sin\theta\cos\theta\right]\r{w}\r{n}[-\cot\theta].\leqno(3.10)
$$ 
Since $f$ vanishes on the small cell, we can restrict ourselves
to the case $\sin\theta\ne0$. We have
$$
\leqalignno{ &f(\r{g})\ll y^{\Re{z_2}-{1\over2}}\left|{x\over
y}-\cot\theta\right|^{-\Re w_1}
\left|{x\over y}+\tan\theta\right|^{-\Re w_2}&(3.11)\cr
&\times|\sin\theta|^{-\Re(w_1+w_3)}|\cos\theta|^{-\Re(w_2+w_4)}
\left|\psi\left({\displaystyle{x\over
y}-\cot\theta\over\displaystyle   {x\over y}+\tan\theta}\right)
\eta(-\cot\theta)\right|. }
$$ 
We claim that if $\Re w_j$ are all bounded then
$$ 
|\sin\theta|^{-w_1-w_3}|\cos\theta|^{-w_2-w_4}\eta(-\cot\theta)
\ll 1,\leqno(3.12)
$$ 
and if moreover $\Re(w_1+w_2)\ge0$ then
$$
\left|{x\over y}-\cot\theta\right|^{- w_1}
\left|{x\over y}+\tan\theta\right|^{-w_2}
\psi\left({\displaystyle{x\over y}-\cot\theta\over\displaystyle
{x\over y}+\tan\theta}\right)\ll
\min\left(1,\left|{x\over
y}\right|^{-\Re(w_1+w_2)}\right).\leqno(3.13)
$$ 
To prove $(3.12)$ we may assume that
$|\cos\theta|\le|\sin\theta|$, and consequentially
$|\cos\theta|\le1/\sqrt{2}\le|\sin\theta|$. Then, by $(3.6)$ the left
side is $\ll|\sin\theta|^{-C-\Re(w_1+w_3)}
|\cos\theta|^{C-\Re(w_2+w_4)}\ll1$. To prove
$(3.13)$ we need to consider the two cases $|x/y|<2$ and
$|x/y|\ge 2$ separately. In the first case we note that
$|x/y-\cot\theta|+|x/y+\tan\theta|
\ge|\cot\theta+\tan\theta|\ge2$. Thus we may assume, for instance,
that $X=|x/y-\cot\theta|\ge1$. We have, with $Y=|x/y+\tan\theta|$,
that either
$Y\ge X\ge1$ or $X\ge Y\ge1$ or $X\ge1\ge Y$. The left side of
$(3.13)$ is, by $(3.2)$,
$$
\ll X^{-\Re w_1}Y^{-\Re w_2}
\min\left(\left({X\over Y}\right)^C,
\left({Y\over X}\right)^C\right)\le1,\leqno(3.14)
$$ 
provided $\Re (w_1+ w_2)\ge0$. In the remaining case the left
side of $(3.13)$ is
$$
\ll \left|{x\over y}\right|^{-\Re(w_1+w_2)}
  X_1^{-\Re w_1}Y_1^{-\Re w_2}\min\left(\left({X_1\over
Y_1}\right)^C,
\left({Y_1\over X_1}\right)^C\right),
\leqno(3.15)
$$ 
where $X_1=|1-(y/x)\cot\theta|$, $Y_1=|1+(y/x)\tan\theta|$. We
may assume, for instance, that $|\tan\theta|\le1$. Then we have
${1\over2}\le Y_1\le{3\over2}$, and thus $(3.15)$ is
$\ll|x/y|^{-\Re(w_1+w_2)}$, which gives $(3.13)$. Summing up, we have
$$ f(\r{g})\ll y^{{\Re z_2}-{1\over2}}
\min\left(1,\left|{x\over
y}\right|^{-\Re(w_1+w_2)}\right),\leqno(3.16)
$$ 
provided $\Re w_j$ are all bounded and $\Re(w_1+w_2)\ge0$. It
should be remarked that this is proved without taking into account
the effect of the factor $\tau(ad)$.
\par The bound $(3.16)$ gives
$$
\sum_{\mu\in\varGamma_\infty}|f(\mu\r{g})|\ll (1+y)y^{{\Re
z_2}-{1\over2}},\quad
\varGamma_\infty=\varGamma\cap N,\leqno(3.17)
$$ 
if $\Re(w_1+w_2)>1$, and it follows that
$$
\e{P}f(\r{g})=\sum_{\gamma\in\varGamma_\infty\!\backslash\varGamma}
\;\sum_{\mu\in\varGamma_\infty}f(\mu\gamma\r{g})\leqno(3.18)
$$ 
is absolutely convergent for any $\r{g}$ if $(3.8)$ holds. In fact
this is the result of comparing $(3.18)$ with the Eisenstein series
$$ 
E_p(\r{g};\nu)=\sum_{\gamma\in\varGamma_\infty\backslash\varGamma}
\phi_p(\gamma\r{g};\,\nu),\leqno(3.19)
$$ 
which converges absolutely for $\Re\nu>{1\over2}$. The assertion
on the convergence of $\e{T}\e{P}f$ is now immediate, and the
formula $(3.9)$ follows readily.
\par 
It remains to prove that $\e{P}f$ is infinitely differentiable  
throughout
$G$.  We may restrict ourselves to the case where none of the
elements   of the matrix $\r{g}$ is equal to $0$, for
$\tau(ad)\not=0$ implies this.   The boundary situation can be
discussed likewise. Computing it explicitly, we see that
$\b{X}f(\r{g})=(d/dt)_{t=0}f\left(\r{g}
\exp(\b{X}t)\right)$, with $\bf X$ as in $(2.9)$ and $\r{g}$ as in  
$(3.1)$, is a linear combination of the five functions
$$
\eqalign{ a&|a|^{-w_1}|b|^{-w_2-1}|c|^{-w_3}
|d|^{-w_4}\psi(ad/(bc))\tau(ad)\eta(d/c),\cr
c&|a|^{-w_1}|b|^{-w_2}|c|^{-w_3}
|d|^{-w_4-1}\psi(ad/(bc))\tau(ad)\eta(d/c),\cr {a\over
b^2c}&|a|^{-w_1}|b|^{-w_2}|c|^{-w_3}
|d|^{-w_4}\psi'(ad/(bc))\tau(ad)\eta(d/c),\cr
ac&|a|^{-w_1}|b|^{-w_2}|c|^{-w_3}
|d|^{-w_4}\psi(ad/(bc))\tau'(ad)\eta(d/c),\cr
&|a|^{-w_1}|b|^{-w_2}|c|^{-w_3}
|d|^{-w_4}\psi(ad/(bc))\tau(ad)\eta'(d/c). }\leqno(3.20)
$$ 
They are majorized by the right side of
$(3.16)$ because of $(3.2)$, $(3.3)$  and $(3.6)$, and $\e{P}\b{X}f$
is absolutely and uniformly convergent throughout $G$, provided
$(3.8)$. That is,
$\b{X}\e{P}f=\e{P}\b{X}f$. Similarly, one may show the same for
$\e{P}\b{Y}f$ and $\e{P}\b{H}f$. This procedure can be repeated
indefinitely, by virtue of the fast decay of the derivatives of
$\psi$, $\tau$, $\eta$. Hence, for any
$\b{u}\in\e{U}$ we have proved that $\b{u}\e{P}f=\e{P}\b{u}f$ in the
pointwise sense, if
$(3.8)$ holds. We end the proof of the lemma.
\medskip 
More precisely, we have, for any fixed $\b{u}$,
$$
\leqalignno{
\b{u}\e{P}f(\r{g})&=\sum_{\mu\in\varGamma_\infty}
\b{u}f(\mu \r{g})
+\sum_{\scr{\gamma\in\varGamma_\infty\backslash\varGamma}\atop
\scr{\gamma\not\in\varGamma_\infty}}
\sum_{\mu\in\varGamma_\infty}\b{u}f(\mu\gamma\r{g})&(3.21)\cr
&=\sum_{\mu\in\varGamma_\infty}\b{u}f(\mu
\r{g})+O\left(y^{{3\over2}-\Re z_2}\right), }
$$ 
for sufficiently large $y$, provided $w$ is as in Lemma 1. By the
Poisson sum formula
$$
\sum_{\mu\in\varGamma_\infty}\b{u}f(\mu
\r{g})=\int_{-\infty}^\infty  
\b{u}f(\r{n}[u]\r{g})du+O\left(y^{{3\over2}-\Re z_2}\right),
\leqno(3.22)
$$ 
for sufficiently large $y$.  To see this, we use that
$\b{u}f(\r{n}[x]
\r{a}[y]\r{k}[\theta])=y^{z_2-{1\over2}}\b{u}f(\r{n}[x/ 
y]\r{k}[\theta])$, and hence
$$
\leqalignno{
\int_{-\infty}^\infty& \b{u}f(\r{n}[x]
\r{a}[y]\r{k}[\theta])\exp(-2\pi imx)dx&(3.23)\cr &=y^{z_2+{1\over2}}
\int_{-\infty}^\infty \b{u}f(\r{n}[x]\r{k}[\theta])
\exp(-2\pi imyx)dx, }
$$ 
which is $\ll (|m|y+1)^{-C}$ via
integration in parts. The relations
$(3.21)$--$(3.23)$ show that $\b{u}\e{P}f$ is not in
$L^2(\varGamma\backslash G)$, in general. Because of this, we
subtract a $\varGamma$-invariant function from $\b{u}\e{P}f$ to have
a square integrable function (an automorphic regularization): We put
$$
\e{P}_0f=\e{P}f-\e{P}\!{}_\infty  f,\leqno(3.24)
$$ 
where
$$
\e{P}_\infty f(\r{g})=
\sum_{\gamma\in\varGamma_\infty\backslash\varGamma}
\int_{-\infty}^\infty f(\r{n}[u]\gamma\r{g})du,\leqno(3.25)
$$ 
the convergence of which follows from $(3.23)$ with $m=0$,
$\b{u}=1$. An examination of the above shows readily that for each
$\b{u}\in\e{U}$
$$
\b{u}\e{P}_0f\in L^2(\varGamma\backslash G),\leqno(3.26)
$$ 
provided $(3.8)$.
\par 
Observe that
$$
\int_{-\infty}^\infty f(\r{n}[u]\r{g})du=\sum_{p=-\infty}^\infty
f_p\phi_p\left(\r{g};z_2\right),\leqno(3.27)
$$ 
with $f_p\ll (|p|+1)^{-C}$. Hence
$$
\e{P}\!{}_\infty  f(\r{g})=\sum_{p=-\infty}^\infty
f_pE_p\left(\r{g};z_2\right).
\leqno(3.28)
$$ 
Also, a computation shows that
$$
\leqalignno{
\int_{-\infty}^\infty f(\r{n}[u]\r{g})&du=
|c|^{w_2-w_3-1}|d|^{w_1-w_4-1}\eta\left({d\over c}\right)&(3.29)\cr
&\times\int_0^\infty {1\over x^{w_1}(x+1)^{w_2}}\psi\left({x\over
x+1}\right)\tau(-x)dx, }
$$ 
with $\r{g}$ as in $(3.1)$. In particular,
$$
\e{P}\!{}_\infty  f(1)=2{\hat\psi}_\tau(0)
\sum_{c=1}^\infty\sum_{\scr{d=1}\atop
\scr{(c,d)=1}}^\infty c^{w_2-w_3-1}d^{w_1-w_4-1}
\eta\left({d\over c}\right),\leqno(3.30)
$$ 
anticipating $(5.20)$.
\medskip 
We are about to invoke the point-wise spectral expansion in
$L^2(\varGamma\backslash G)$, which is to be applied to $\e{P}_0f$.  
Then we need to quote the Fourier expansion of the Eisenstein series:
$$
\leqalignno{
E_p(\r{g};\nu)&=\phi_p(\r{g};\nu)+d_p(\nu)\phi_p(\r{g};-\nu)
&(3.31)\cr
&+{1\over\zeta(1+2\nu)}
\sum_{\scr{n=-\infty}\atop\scr{n\ne0}}^\infty
|n|^{-{1\over2}-\nu}\sigma_{2\nu}(|n|)
\e{A}^{\sgn(n)}\phi_p(\r{a}[|n|]\r{g};\nu), }
$$ 
provided the right side is finite, where
$$ 
d_p(\nu)={(-1)^p\pi\Gamma(2\nu)\zeta(2\nu)\over
2^{2\nu-1}\zeta(2\nu+1)
\Gamma({1\over2}+\nu+p)\Gamma({1\over2}+\nu-p)}.\leqno(3.32)
$$ 
We have also the functional equation
$$ 
E_p(\r{g};\nu)=d_p(\nu)E_p(\r{g};-\nu),\quad
d_p(\nu)d_p(-\nu)=1.\leqno(3.33)
$$ 
The proof of $(3.31)$--$(3.33)$ is practically the same as that of
Lemma 1.2 of [19].
\medskip 
Making the assertions $(2.10)$, $(2.11)$, and $(2.13)$ precise, we
have
\medskip
\noindent 
{\bf Lemma 2.}\quad {\it Let $\varpi_V$ be the orthogonal
projection to $V$, and $\varpi_E$ to ${}^e\!L^2(\varGamma\backslash
G)$. Let
$h$ be a vector such that $\b{u}h\in L^2(\varGamma\backslash G)$ for
any $\b{u}\in\e{U}$. Then the spectral decomposition
$$ 
h(\r{g})={3\over\pi}
\langle h,1\rangle_{\varGamma\backslash G}+
\sum_V\varpi_Vh(\r{g})+\varpi_Eh(\r{g})
\leqno(3.34)
$$ 
converges absolutely for each $\r{g}\in G$. Similarly
$$
\varpi_Vh(\r{g})=\sum_{p=-\infty}^\infty\langle  
h,\varphi_p\rangle_{\varGamma
\backslash G}\varphi_p(\r{g}),\leqno(3.35)
$$ where $\varphi_p$ are as above together with an obvious convention
for $V$ in the discrete series. Also
$$
\varpi_Eh(\r{g})=
\sum_{p=-\infty}^\infty\int_{(0)}e_p(h;\nu)E_p(\r{g};\nu)
{d\nu\over4\pi i},\leqno(3.36)
$$ 
with
$$ e_p(h;\nu)=\int_{\varGamma\backslash G}h(\r{g})
\overline{E_p(\r{g};\nu)}d\r{g}.\leqno(3.37)
$$ 
}
\par
\noindent 
{\it Proof.\/} This assertion is taken from Section 1.2 of
[8],  which is based on [9]. Other approaches are possible. See the
remark at the end of this section.
\medskip Hence, we have, by  $(3.26)$, a pointwise spectral
decomposition of
$\e{P}_0f$.  We may put the result as
$$
\e{P}f(\r{g})=\e{P}\!{}_\infty  f(\r{g})+
\sum_V\varpi_V\e{P}_0f(\r{g})+\varpi_E\e{P}_0f(\r{g}),
\leqno(3.38)
$$ 
with
$$
\varpi_E\e{P}_0f(\r{g})=
\sum_{p=-\infty}^\infty\int_{(0)}e_p^{(1)}(\e{P}_0
f;\nu)E_p(\r{g};\nu) {d\nu\over4\pi i}.\leqno(3.39)
$$ 
Here, $e_p^{(1)}$ is the part of $e_p$ corresponding to the third
term on the right of
$(3.31)$. The identity $(3.38)$ depends on the fact that $\langle
\e{P}_0f,1\rangle_{\varGamma
\backslash G}=0$ and
$e_p(\e{P}_0f;\nu)=e_p^{(1)}(\e{P}_0f;\nu)$ for all
$p$. Both follow readily from the definition $(3.24)$ of $\e{P}_0f$.
\par The last sum can be taken inside the integral because of
absolute   convergence (see the remark given after $(5.33)$ below).
Applying the Hecke operator to $(3.38)$ and
$(3.39)$ termwise, and invoking $(2.26)$ and $(3.28)$, we get, on  
$(3.8)$,
$$
\e{T}\e{P}f(\r{g})=\e{T}\e{P}\!{}_\infty  f(\r{g})+
\sum_V\e{T}\varpi_V\e{P}_0f(\r{g})+\e{T}(\varpi_E^{(0)}+
\varpi_E^{(1)})\e{P}_0f(\r{g}),\leqno(3.40)
$$ 
where
$$
\eqalign{
\e{T}\e{P}_\infty f(\r{g})&=
\zeta(w_1+w_2-1)\zeta(w_3+w_4)\e{P}\!{}_\infty  f(\r{g}),\cr
\e{T}\varpi_V\e{P}_0f(\r{g})&=
H_V\left(z_1\right)\varpi_V\e{P}_0f(\r{g}),\cr
\e{T}\varpi^{(j)}_E\e{P}_0f(\r{g})&=
\int_{(0)}\zeta\left(z_1+\nu\right)
\zeta\left(z_1-\nu\right)
\e{E}^{(j)}(\e{P}_0f;\r{g},\nu){d\nu\over4\pi i}, }\leqno(3.41)
$$ 
with $z_1$ as in $(2.3)$. Here we have used
$T_nE_p(\r{g};\nu)=n^{-\nu}\sigma_{2\nu}(n)E_p(\r{g};\nu)$, and have
put
$$
\e{E}^{(j)}(\e{P}_0f;\r{g},\nu)=\sum_{p=-\infty}^\infty
e_p^{(1)}(\e{P}_0f;\nu)E_p^{(j)}(\r{g};\nu),\leqno(3.42)
$$ 
where $E^{(0)}_p$ is the sum of the first two terms on the right
of
$(3.31)$ and $E^{(1)}_p$ the rest. Observe that the three Hecke
series in $(3.41)$, i.e., $H_V(z_1)$ and its analogues, are all
absolutely bounded under $(3.8)$.
\par In view of $(3.9)$, we use
$(3.40)$ when $\r{g}$ is the unit element. Namely, our original
problem has been reduced to computing the quantities
$\varpi_V\e{P}_0f(1)$ and $\e{E}^{(j)}(\e{P}_0f;1,\nu)$. Note that we
have already $(3.30)$.
\medskip
\noindent 
{\bf Remark.} The spectral decomposition in
$L^2(\varGamma\backslash G)$ can be derived, via the Fourier
expansion with respect to the right action of $K$, from that in
$L^2(\varGamma\backslash \B{H})$, thus for instance, from a minor
extension of [19, Chapter 1]. Naturally, the pointwise convergence
in $(3.38)$ is crucial for our purpose. Thus, it should be stressed
that there is a way, based on explicit estimation, to achieve the
same without recourse to Lemma 2 but rather starting with the
convergence in $L^2(\varGamma\backslash G)$. The necessary
uniform estimate is in fact provided by the discussion of either
Section 5 or 6. This indicates that we may relax
considerably the stringent decay condition on derivatives of
$\psi$,  $\tau$, $\eta$ imposed above. 
\vskip 1cm
\centerline
{\bf 4. Big cell}
\bigskip The unfolding argument reduces our task further to an
application of the harmonic analysis in the big cell of the Bruhat
decomposition, as is to be seen in the next section.  Hence we
collect here fundamentals in this context, which may be termed the
Kirillov scheme collectively.
\medskip We first extend $(2.15)$ by
$$
\e{A}^\delta\phi=\sum_{p=-\infty}^\infty c_p\e{A}^\delta\phi_p,
\quad \phi=\sum_{p=-\infty}^\infty c_p\phi_p,\leqno(4.1)
$$ 
where $\phi_p(\r{g})=\phi_p(\r{g};\nu)$, and $\phi$ is smooth,
i.e., $c_p$ decays faster than any negative power of $|p|$ as
$p$ tends to infinity. Note that  we occasionally omit to mention
$\nu$.  We shall show that $(4.1)$ exists for any $\nu$. For this
and other purposes, the following estimates will be useful; bounds
up to
$(4.5)$ are all uniform for $p$ and $|\Re\nu|<{1\over2}$.
\par The first line of $(2.16)$ gives
$$
\leqalignno{
\e{A}^\delta\phi_p(\r{a}[y])&=
\e{A}^\delta\phi_0(\r{a}[y])+y^{{1\over2}-\nu}\int_{-\infty}^\infty
{\exp(2\pi iy\xi)\over (\xi^2+1)^{{1\over2}+\nu}}
\left(\left({\xi+i\over\xi-i}\right)^{\delta p}-1\right)
d\xi&(4.2)\cr
&={2\pi^{{1\over2}+\nu}\over\Gamma({1\over2}+\nu)}y^{1\over2}
K_\nu(2\pi y) +O\left(y^{{1\over2}-\Re\nu}(|p|+1)\right). }
$$ 
By the power series expansion for
$K_\nu$, we get, as $y\to0^+$,
$$
\e{A}^\delta\phi_p(\r{a}[y])\ll
(|p|+|\nu|+1)y^{{1\over2}-|\Re\nu|-\varepsilon}.\leqno(4.3)
$$ 
On the other hand, we have, by integration in parts,
$$
\e{A}^\delta\phi_p(\r{a}[y])= {y^{-{1\over2}-\nu}\over2\pi i}
\int_{-\infty}^\infty{((1+2\nu)\xi+2\delta pi)\exp(2\pi iy\xi)
\over(\xi^2+1)^{{3\over2}+\nu}}
\left({\xi+i\over\xi-i}\right)^{\delta p} d\xi.\leqno(4.4)
$$ 
Shifting the contour to $\Im\xi=(|\nu|+|p|+1)^{-1}$, we see that
$$
\e{A}^\delta\phi_p(\r{a}[y])\ll
(|p|+|\nu|+1)y^{-{1\over2}-\Re\nu}\exp\left(-{y\over  
|\nu|+|p|+1}\right).
\leqno(4.5)
$$ 
Repeating integration in parts in $(4.4)$, we find that $(4.1)$
converges in any fixed vertical strip of the $\nu$-plane. Note that
$$
\e{A}^\delta\phi(\r{g})=\int_{-\infty}^\infty\exp(-2\pi\delta
ix)\phi(\r{w}\r{n}[x]\r{g})dx,\leqno(4.6)
$$ 
for those $\nu$ in the domain where the integral converges
uniformly. In fact the equality holds at least for $\Re\nu>0$, and
the assertion follows by analytic continuation.
\medskip We then define the Kirillov operator $\e{K}$ by
$$
\e{K}\phi(u)=\e{A}^{\sgn(u)}\phi(\r{a}[|u|]),\quad u\in\B{R}^\times.
\leqno(4.7)
$$ 
This concept will play a crucial r\^ole in our argument, via the
following three lemmas:
\medskip
\noindent {\bf Lemma 3.}\quad{\it Let $\phi$ be smooth as in
$(4.1)$. We have, with the right translation $R$,
$$
\e{K}R_{\r{n}[x]}\phi(u)=\exp(2\pi iux)\e{K}\phi(u),\quad
\e{K}R_{\r{a}[y]}\phi(u)=\e{K}\phi(uy).\leqno(4.8)
$$ 
Also, if $|\Re\nu|<{1\over2}$, then
$$
\e{K}R_\r{w}\phi(u)=\int_{\B{R}^\times}j_\nu(u\lambda)\e{K}
\phi(\lambda)d^\times\!\lambda,\leqno(4.9)
$$ 
where $j_\nu$ is defined by $(2.33)$. }
\medskip
\noindent 
{\it Proof.\/}  The explicit description $(4.9)$ of the action 
of the Weyl element is probably due
to N.Ja.\ Vilenkin (see Section 7 of [21, Chapter VII] as
well as the   formula (17) of [22, p.\ 454]). A rigorous proof can
be found in Theorem 2 of [20], which is developed in the context of
automorphy but in fact asserts the above. It is shown there that the
function
$$
\Gamma_p(s)=\int_0^\infty
\e{A}^+\phi_p(\r{a}[y])y^{s-{3\over2}}dy\leqno(4.10)
$$ 
continues meromorphically to $\B{C}$, and satisfies the
Jacquet--Langlands local functional equation
$$
\leqalignno{ (-1)^p\Gamma_p(s) =&2^{1-2s}\pi^{-2s}\Gamma(s+\nu)
\Gamma(s-\nu)&(4.11) \cr &\times\left(\cos\pi
s\,\Gamma_p(1-s)+\cos\pi\nu\,
\Gamma_{-p}(1-s)\right). }
$$ 
The Mellin inversion of this coupled with $(2.39)$--$(2.40)$
gives  
$(4.9)$ for
$\phi=\phi_p$.  A combination of $(2.45)$, $(4.3)$, and $(4.5)$
yields   the necessary analytic continuation in $\nu$, and the
extension to smooth  $\phi$.
\medskip
\noindent {\bf Lemma 4.}\quad{\it Let $\nu\in i\B{R}$, and introduce
the Hilbert space
$$ 
U_\nu=\overline{\bigoplus_{p=-\infty}^\infty
\B{C}\phi_p},\quad
\phi_p(\r{g})=\phi_p(\r{g};\nu),
\leqno(4.12)
$$ 
equipped with the norm
$$
\Vert\phi\Vert_{U_\nu}=\sqrt{\sum_{p=-\infty}^\infty
|c_p|^2},\quad\phi= \sum_{p=-\infty}^\infty c_p\phi_p.\leqno(4.13)
$$ 
Then the operator $\e{K}$ is a unitary map from $U_\nu$ onto
$L^2(\B{R}^\times,
\pi^{-1}d^\times)$.
\/}
\medskip
\noindent 
{\it Proof.\/} This seems to stem from A.A. Kirillov [15].
A proof of the unitaricity is given in Theorem 1 of [20], though
disguised in the context of automorphy. It depends on the following
integral formula: For any
$\alpha,\beta\in\B{C}$ and
$|\Re\nu|<{1\over2}$
$$
\leqalignno{
\qquad\int_0^\infty& W_{\alpha,\nu}(u) W_{\beta,\nu}(u){du\over u}=
{\pi\over(\alpha-\beta)\sin(2\pi\nu)}&(4.14)\cr
&\times\left[{1\over\Gamma({1\over2}-\alpha+\nu)\Gamma({1\over2}
-\beta-\nu)}-{1\over\Gamma({1\over2}-\alpha-\nu)\Gamma({1\over2}
-\beta+\nu)}\right] }
$$ 
(see the formula 7.611(3) of [11]). The proof in [20] of this
employs   the Whittaker differential equation.  Here we shall show
the surjectivity   of the map. Thus, let us assume that $\nu\in
i\B{R}$, and that a smooth   function
$k$, compactly supported on $\B{R}^\times$, is orthogonal to all
$\e{K}\phi_p$. Multiply $(4.4)$ by $k$ and integrate, change the
order of integration, and undo the integration in parts with respect
to the outer integral. We have
$$
\leqalignno{ 0&=\int_{\B{R}^\times}k(u)\overline{\e{K}\phi_p(u)}
d^\times\!u&(4.15)\cr &=\int_{-\infty}^\infty
{1\over(\xi^2+1)^{{1\over2}+\nu}}
\left({\xi+i\over
\xi-i}\right)^p\int_{-\infty}^\infty k(u) |u|^{-{1\over2}+\nu}
\exp(-2\pi i u\xi)du\,d\xi. }
$$ 
Observe that the system
$\left\{((\xi+i)/(\xi-i))^p:\,p\in\B{Z}\right\}$ is complete
orthonormal in the space  
$L^2\left(\B{R},(\pi(\xi^2+1))^{-1}d\xi\right)$. Hence the Fourier
transform of $k(u)|u|^{-{1\over2}+\nu}$ vanishes identically, whence
the assertion.
\medskip 
Next, we consider the complementary series or the
situation with
$-{1\over2}<\nu<{1\over2}$. This is included here only for the sake
of completeness; such a representation of $G$ does not occur in
$L^2(\varGamma\backslash G)$.  Obviously, Lemma 2 remains   valid.
The definition $(4.12)$ is the same, but $(4.13)$ has to be replaced
by the   norm
$$
\sqrt{\pi^{2\nu}\sum_{p=-\infty}^\infty{\Gamma(p+{1\over2}-\nu)
\over\Gamma(p+{1\over2}+\nu)} |c_p|^2}\,.\leqno(4.16)
$$ 
With this, the above proof extends readily, and Lemma 4 holds for
these $\nu$ as well.
\medskip 
On the other hand, in dealing with the holomorphic discrete
series,  
$(4.12)$ needs to be replaced by the Hilbert space
$$ D_\ell=\overline{\bigoplus_{p=\ell}^\infty
\B{C}\phi_p},\quad \phi_p(\r{g})=
\phi_p\left(\r{g};\ell-{1\over2}\right),\quad  
1\le\ell\in\B{Z}\,,\leqno(4.17)
$$ 
equipped with the norm
$$
\Vert\phi\Vert_{D_\ell}=\sqrt{\pi^{2\ell-1}\sum_{p=\ell}^\infty
{\Gamma(p-\ell+1)\over\Gamma(p+\ell)}|c_p|^2},
\quad \phi=\sum_{p=\ell}^\infty c_p\phi_p.\leqno(4.18)
$$ 
Since $\e{A}^-$ annihilates $D_\ell$, we are concerned with
$\e{A}^+$ only. The expression $(2.16)$, $\delta=+$, holds without
changes. With this, the operator $\e{K}$ is defined as before.
\medskip
\noindent 
{\bf Lemma 5.} {\it The operator $\e{K}$ is a unitary map
from $D_\ell$ onto $L^2((0,\infty),\pi^{-1}d^\times)$. Also, for any
smooth $\phi\in D_\ell$, we have $(4.9)$ with $\nu=\ell-{1\over2}$.
The analogue for the antiholomorphic discrete series is obtained by  
applying the involution $\omega$. 
}
\medskip
\noindent 
{\it Proof.\/} The third assertion is immediate. As to the
unitaricity of $\e{K}$, it is proved with a minor change of the
above argument. In fact, the Whittaker function
$W_{p,\ell-{1\over2}}(u)$ ($p\ge l$) is a product of
$\exp(-u/2)u^\ell$ and a polynomial on $u$ of degree $p-\ell$, as
$(2.16)$ implies. Thus the proof of $(4.14)$ in [20] can be carried
out also for the product
$W_{p,\ell-{1\over2}}(u)W_{q,\ell-{1\over2}}(u)$ with integers
$p,\,q$, although the condition on $\Re\nu$ is violated. The result
is equal to the limit of $(4.14)$ as $(\alpha,\beta,\nu)$ tends to
$\left(p,q,\ell-{1\over2}\right)$.  As to the surjectivity, we
argue as follows: Let
$k$ be smooth and compactly supported on
$(0,\infty)$. If $k$ is orthogonal to all
$\e{K}\phi_p$, $\ell\le p$, then we have, by the remark just made on
$W_{p,\ell-{1\over2}}(u)$,
$$
\int_0^\infty k(u)\exp(-2\pi u)u^p{du\over u}=0,\quad \ell\le p.
\leqno(4.19)
$$ 
This implies that the Fourier transform of
$k(u)\exp(-2\pi u)u^{\ell-1}$ vanishes identically; in fact it
suffices to expand the additive character into a power series and
integrate termwise. Hence $k\equiv0$. The counterpart of $(4.9)$,
with $\phi=\phi_p$, can be proved in much the same way as before. Its
extension to any smooth $\phi$ is immediate with $(2.46)$ and
$$
\e{K}\phi_p(u)=\e{A}^+\phi_p(\r{a}[u])
\ll\min(u,|p|+1)u^{-\ell},\quad u>0,\leqno(4.20)
$$ 
uniformly in $p$. This comes from $(2.16)$ and $(4.4)$.
\medskip
\noindent 
{\bf Remark.} The identity $(4.9)$ is crucial for our
purpose. In a context related to ours, this is given in Theorem 4.1
of [8] without proof nor attribution. The formula seems to have
appeared in print for the first time in [21], thus our partial
attribution above, but the argument there lacks an adequate
discussion on the convergence issue; the same can be said about
[22]. A rigorous proof is given in [20], which is outlined
above; the argument depends on the Mellin transform, and thus is
different from those in  [21], [22]. On the other hand, the recent
article [2] provides an independent proof along the line of [21].
The present authors thank M. Baruch for this information. We stress
that Lemmas 3 and 4 extend to
$\r{PSL}_2(\B{C})$; see [6, Part XIII], which itself is an extension
of [20]. The corresponding surjectivity assertion is not discussed
in [6, Part XIII], but the above argument starting with $(4.15)$
should extend to
$\r{PSL}_2(\B{C})$ on the basis of the discussion in [7, Section 5].
\medskip
\vskip 1cm
\centerline{\bf 5. Projections}
\bigskip 
We are now ready to deal with the task encountered at the
end of Section 3 or the explicit computation of
$\varpi_V\e{P}_0f(1)$ and
$\e{E}^{(j)}(\e{P}_0f;1,\nu)$ in terms of $\psi$, $\tau$,
$\eta$. The basic implement is the Kirillov scheme. The condition  
$(3.8)$ is imposed throughout the present section.
\medskip Let us first consider $V$ in the unitary principal series,
so that
$\nu_V\in i\B{R}$. Let $\varphi_p$ be as in $(2.14)$. Since
Eisenstein series are orthogonal to any cuspidal element, we have
$$
\leqalignno{
\langle\e{P}_0f,\varphi_p\rangle_{\varGamma\backslash G}&=
\langle\e{P}f,\varphi_p\rangle_{\varGamma\backslash G}&(5.1)\cr
&=\int_G f(\r{g})\overline{\varphi_p(\r{g})}d\r{g}. }
$$ 
The unfolding procedure in the second line is justified by
$(3.16)$ and the exponential decay of $\varphi_p$ as $y\to\infty$.
The latter follows from $(2.30)$ and $(4.5)$. We have
$$
\leqalignno{
\langle\e{P}_0f,\varphi_p\rangle_{\varGamma\backslash G}&=
\sum_{\scr{n=-\infty}\atop\scr{n\ne0}}^\infty
{\overline{\varrho_V(n)}\over\sqrt{|n|}}\int_G
f(\r{g})\overline{\e{A}^{\sgn(n)}\phi_p(\r{a}[|n|]\r{g}})d\r{g}&(5.2)\cr
&=H_V\left(z_2\right)
\left(\overline{\varrho_V(1)}\Phi_p^+
+\overline{\varrho_V(-1)}\Phi^-_p\right)f(\nu_V), }
$$ 
where $(2.27)$ and $(3.7)$ have been used, and
$$
\Phi_p^\delta f(\nu)=\int_G
f(\r{g})\overline{\e{A}^\delta\phi_p(\r{g};\nu)}d\r{g}.\leqno(5.3)
$$ 
The absolute convergence that is necessary to have the first line
of $(5.2)$ follows from that of $(5.3)$, which in turn results from
$(2.30)$ and $(4.5)$. By $(3.35)$ or gathering together the
projections of $\e{P}_0f$ to $V_p$, we have now
$$
\leqalignno{ &\varpi_V\e{P}_0f(\r{g})={1\over2}
\left(|\varrho_V(1)|^2+|\varrho_V(-1)|^2\right)
H_V\left(z_2\right)&(5.4)\cr
&\times\sum_{n=1}^\infty{t_V(n)\over\sqrt{n}}
\left(\e{B}^{(+,+)}+\e{B}^{(-,-)}+\epsilon_V\e{B}^{(+,-)}
+\epsilon_V\e{B}^{(-,+)}\right)f(\r{a}[n]\r{g};\,\nu_V), }
$$ 
with
$$
\leqalignno{
\e{B}^{(\delta_1,\delta_2)}
f(\r{g};\nu)&=\sum_{p=-\infty}^\infty\Phi_p^{\delta_1}f(\nu)\,
\e{A}^{\delta_2}\phi_p(\r{g};\nu)&(5.5)\cr &=\exp(2\pi i\delta_2
x)\sum_{p=-\infty}^\infty\Phi_p^{\delta_1}
f(\nu)\e{A}^{\delta_2}\phi_p(\r{a}[y])\exp(2ip\theta). }
$$ 
We shall prove that the right side of $(5.4)$ converges
absolutely to a continuous function in $V$.
\medskip 
To this end, we show the bound
$$
\Phi_p^\delta f(\nu)\ll (|p|+|\nu|+1)^{-C},\quad
|\Re\nu|<{1\over2}.\leqno(5.6)
$$ 
A combination of $(4.5)$ and $(5.6)$ yields
$$
\e{B}^{(\delta_1,\delta_2)}f(\r{g};\nu)\ll
y^{{1\over2}-|\Re\nu|-\varepsilon} ((y+1)(|\nu|+1))^{-C},\leqno(5.7)
$$ 
in the same region of $\nu$, whence the   above claim on $(5.4)$.
To prove $(5.6)$, observe, as in the proof of Lemma 1, that the
function $\b{u}f$ is bounded by the right side of $(3.16)$, for any
given $\b{u}\in\e{U}$. Thus, the second line of $(2.16)$ and $(4.5)$
give
$$
\leqalignno{
\Phi_p^\delta\b{u}f&\ll\int_0^\infty y^{\Re  
z_2-{3\over2}}\left|\e{A}^\delta\phi_p(\r{a}[y])\right|dy&(5.8)\cr
&\ll(|p|+|\nu|+1)^{\Re z_2-\Re\nu}. }
$$ 
Since $\e{A}^\delta$ is an intertwining operator with respect to
the action of the elements of $\e{U}$, we have, for any positive
integer $q$,
$$
\Phi_p^\delta(\Omega+ i\partial_\theta^2)^qf=
\left(\nu^2-{1\over4}-4ip^2\right)^q\Phi_p^\delta f,\leqno(5.9)
$$ 
by integration in parts, which can be justified with
$(5.8)$. This obviously gives $(5.6)$.
\medskip 
We may now look at $\e{B}^{(\delta_1,\delta_2)}
f(\r{a}[y];\nu)$ closer with the Kirillov scheme: We assume that
$\nu\in i\B{R}$. We have
$$
\e{B}^{(\delta_1,\delta_2)} f(\r{a}[y];\nu)=\sum_{p=-\infty}^\infty
\Phi_p^{\delta_1}f(\nu)\,\e{K}\phi_p(\delta_2y)
=\e{K}\e{L}^{\delta_1}f(\delta_2 y),
\leqno(5.10)
$$ 
where
$$
\e{L}^\delta
f=\sum_{p=-\infty}^\infty\Phi_p^{\delta}f(\nu)\,\phi_p\leqno(5.11)
$$ 
is a smooth element in $U_\nu$, because of $(5.6)$. The
unitaricity assertion in Lemma 4 gives
$$
\Phi_p^{\delta}f=\langle \e{L}^\delta f,\phi_p\rangle_{U_\nu}
={1\over\pi}\int_{\B{R}^\times}
\e{K}\e{L}^\delta f(u)\overline{\e{K}\phi_p(u)} d^\times\!
u.\leqno(5.12)
$$ 
This and the surjectivity assertion there imply that if one can  
transform
$(5.3)$ into
$$
\Phi_p^{\delta}f={1\over\pi}\int_{\B{R}^\times}Y^\delta(u)
\overline{\e{K}\phi_p(u)}d^\times\! u\leqno(5.13)
$$ 
then it should follow that
$$
\e{B}^{(\delta_1,\delta_2)} f(\r{a}[y])=Y^{\delta_1}(\delta_2
y).\leqno(5.14)
$$ 
Note that we have used implicitly a simple continuity argument,
which will be made explicit at $(5.22)$.
\medskip We may write $(5.3)$ as
$$
\Phi_p^{\delta}f =\int_0^\infty
u^{z_2-{3\over2}}\int_{N\r{w}N}f(\r{g})
\overline{R_\r{g}\e{A}^\delta\phi_p(\r{a}[u])}d\dot\r{g}du.
\leqno(5.15)
$$ 
Here $\r{g}=\r{n}[x_1]\r{w}\r{n}[x_2]$ and
$d\dot\r{g}=\pi^{-1}dx_1dx_2$; the formula $(3.10)$ gives the
Jacobian   for this change of variables. We observe
$$
\leqalignno{ R_\r{g}\e{A}^\delta\phi_p(\r{a}[u])&=R_{\r{wn}[x_2]}
\e{A}^\delta\phi_p(\r{n}[x_1u]\r{a}[u])=\exp(2\pi i\delta x_1u)
R_{\r{wn}[x_2]}\e{A}^\delta\phi_p(\r{a}[u])&(5.16)\cr &=\exp(2\pi
i\delta x_1u)\e{A}^\delta R_\r{w}R_{\r{n}[x_2]}
\phi_p(\r{a}[u]). }
$$ 
By Lemma 3
$$
\leqalignno{
\e{A}^\delta R_\r{w}R_{\r{n}[x_2]}\phi_p(\r{a}[u])
&=\e{K}R_\r{w}R_{\r{n}[x_2]}\phi_p(\delta u)=\int_{\B{R}^\times}
j_\nu(\delta
u\lambda)\e{K}R_{\r{n}[x_2]}\phi_p(\lambda)d^\times\!\lambda&(5.17)\cr
&=\int_{\B{R}^\times} \exp(2\pi ix_2\lambda)j_\nu(\delta
u\lambda)\e{K}\phi_p(\lambda)d^\times\!\lambda. }
$$ 
Thus
$$
\leqalignno{
\Phi_p^{\delta}f=&{1\over\pi}\int_0^\infty  
u^{z_2-{3\over2}}\int_{\B{R}^2}
f(\r{n}[x_1]\r{w}\r{n}[x_2])\exp(-2\pi i\delta x_1u)&(5.18)\cr
&\times\int_{\B{R}^\times} \exp(-2\pi ix_2\lambda)j_\nu(\delta
u\lambda)\overline{\e{K}\phi_p(\lambda)}d^\times\!\lambda dx_1dx_2du,
}
$$ 
where we have used that $j_\nu=j_{-\nu}$. Applying change of
variables
$x_1\to \displaystyle{x_1\over x_2}$,  $x_2\to -x_2$, we have
$$
\leqalignno{
\Phi_p^\delta f=&{1\over\pi}\int_0^\infty
u^{z_2-{3\over2}}\int_{\B{R}^\times}
|x_2|^{w_1-w_4}{\hat\psi}_\tau\left(\delta{ u\over x_2}\right)
\eta(x_2)&(5.19)
\cr &\times\int_{\B{R}^\times}\exp\left(2\pi i{\lambda
x_2}\right)j_\nu(\delta u\lambda)\overline{\e{K}\phi_p(\lambda)}
d^\times\!\lambda\, d^\times\!x_2\, du, }
$$ 
with
$$ {\hat\psi}_\tau(u)=\int_0^\infty {1\over x_1^{w_1}(x_1+1)^{w_2}}
\psi\left({x_1\over x_1+1}\right)
\tau(-x_1)\exp(-2\pi iux_1)dx_1.\leqno(5.20)
$$ 
Here $(3.3)$ and $(3.6)$ have been used. The triple integral in
$(5.19)$ converges absolutely. In fact, a multiple application of
integration in parts gives, for each $l$,
$$
\left({d\over du}\right)^l
\hat{\psi}_\tau(u)\ll(|u|+1)^{-C},\leqno(5.21)
$$ 
because of $(3.3)$. A combination of $(2.45)$, $(4.3)$, $(4.5)$,
and $(5.21)$ yields that the integral whose integrand is the
absolute value of that in $(5.19)$ is $\ll|p|+1$, with the implied
constant depending on
$\nu$. Hence we have, for any smooth $\phi\in U_\nu$,
$$
\leqalignno{
\langle\e{L}^\delta f,\phi\rangle_{U_\nu}&={1\over\pi}
\int_{\B{R}^\times}\Bigg\{
\int_0^\infty  u^{z_2-{3\over2}}\int_{\B{R}^\times}
|x_2|^{w_1-w_4}&(5.22)\cr &\times j_\nu(\delta
u\lambda){\hat\psi}_\tau\left(\delta {u\over x_2}\right)\eta(x_2)
\exp\left(2\pi i{\lambda  x_2}\right)\, d^\times\!x_2\, du\Bigg\}
\overline{\e{K}\phi(\lambda)}d^\times\!\lambda. }
$$ 
Via $(5.10)$, Lemma 4 now gives rise to
$$
\e{B}^{(\delta_1,\delta_2)}f(\r{a}[y];\nu)=
\e{B}^{\delta_1\delta_2}f(\r{a}[y];\nu),\leqno(5.23)
$$ 
with
$$
\leqalignno{
\e{B}^\delta &f(\r{a}[y];\nu) =\int_0^\infty u^{z_2-{3\over2}}
j_\nu(\delta uy)&(5.24)\cr &\times\int_{\B{R}^\times} |x_2|^{w_1-w_4}
{\hat\psi}_\tau\left(\delta
{u\over x_2}\right)\eta(x_2)\exp\left( 2\pi i{y
x_2}\right)d^\times\!x_2\,du. }
$$ 
Inserting this into $(5.4)$, we obtain, via $(3.41)$,
\medskip
\noindent {\bf Lemma 6.} {\it If $V$ is in the unitary principal
series, then
$$
\leqalignno{
\e{T}\varpi_V\e{P}_0f(1)
=&\left(|\varrho_V(1)|^2+|\varrho_V(-1)|^2\right)H_V(z_1)
H_V\left(z_2\right) &(5.25)\cr
&\times\sum_{n=1}^\infty{t_V(n)\over\sqrt{n}}
\left(\e{B}^++\epsilon_V\e{B}^-\right) f(\r{a}[n];\,\nu_V), 
}
$$ 
provided $(3.8)$. }
\medskip 
Next, we treat the discrete series. We assume that $V$ is
in the   holomorphic discrete series, having the complete
orthonormal system
$\{\varphi_p:\,p\ge\ell\}$ with $\varphi_p$ given in $(2.21)$. The
relation $(5.1)$ extends as it is; the unfolding procedure depends
on the observation on the Whittaker function $W_{p,\ell-{1\over2}}$
made in the proof of Lemma 5. Then,
$(5.2)$, with an obvious interpretation of
$(5.3)$, is replaced by
$$
\langle \e{P}f,\varphi_p\rangle_{\varGamma\backslash G}=
\pi^{{1\over2}-\ell}\overline{\varrho_V(1)}H_V\left(z_2\right)
\left({\Gamma(p+\ell)\over\Gamma(p-\ell+1)}\right)^{1\over2}
\Phi_p^+f\left(\ell-{1\over2}\right),\leqno(5.26)
$$ 
and $(5.4)$ by
$$
\varpi_V\e{P}_0f(\r{g})=|\varrho_V(1)|^2H_V\left(z_2\right)
\sum_{n=1}^\infty{t_V(n)\over\sqrt{n}}
\e{B}f\left(\r{a}[n]\r{g};\ell-{1\over2}\right).\leqno(5.27)
$$ 
Here
$$
\e{B}f\left(\r{g};\ell-{1\over2}\right)=\pi^{1-2\ell}
\sum_{p=\ell}^\infty
{\Gamma(p+\ell)\over\Gamma(p-\ell+1)}\Phi_p^+f\left(\ell- 
{1\over2}\right)
\e{A}^+\phi_p\left(\r{g};\ell-{1\over2}\right),\leqno(5.28)
$$ 
which replaces $(5.5)$. The $\e{B}f$ exists  as a continuous
function in $V$. On noting $(4.20)$, this follows from
$\Phi_p^+f(\ell-{1\over2})\ll (|p|+1)^{-C}$. To get the latter, we
observe that for any
$\b{u}\in\e{U}$
$$
\langle \e{P}\b{u}f,\varphi_p\rangle_{\varGamma\backslash G}=\langle
\b{u}\e{P}f,\varphi_p\rangle_{\varGamma\backslash G}=\pm\langle
\e{P}f,\overline\b{u}\varphi_p\rangle_{\varGamma\backslash
G}.\leqno(5.29)
$$ 
Set $\b{u}=\partial_\theta^q$, with a positive integer $q$; and
use $(5.26)$ on the right side and  $|\langle
\e{P}\b{u}f,\varphi_p\rangle_{\varGamma\backslash G}|\le\Vert
\e{P}_0\b{u}f\Vert_{\varGamma\backslash G}$ on the left, which
confirms the claim. We have actually proved that $\varpi_V\e{P}_0f$
exists as a continuous function in $V$.
\medskip 
We prove an extension of $(5.24)$. This is now easy: We
put, in place of $\e{L}^\delta f$,
$$
\e{L}f=\pi^{1-2\ell}\sum_{p=\ell}^\infty
{\Gamma(p+\ell)\over\Gamma(p-\ell+1)}\Phi_p^+f
\left(\ell-{1\over2}\right)
\phi_p,\leqno(5.30)
$$ 
which is a smooth element in $D_\ell$. Then, we can proceed much
like $(5.10)$--$(5.22)$, relying on Lemma 5 and $(2.46)$, $(4.20)$.
Thus, we have
$$
\e{B}f(\r{a}[y])=\e{B}^+f\left(\r{a}[y];\ell-{1\over2}\right),
\leqno(5.31)
$$ 
with an extended use of notation.
\medskip Hence, taking into account the antiholomorphic discrete
series as well,   we have
\medskip
\noindent 
{\bf Lemma 7.} {\it If $V$ is in the discrete series, then
$$
\e{T}\varpi_V\e{P}_0f(1)=\left(|\varrho_V(1)|^2+|
\varrho_V(-1)|^2\right)
H_V(z_1)H_V\left(z_2\right)
\sum_{n=1}^\infty{t_V(n)\over\sqrt{n}}
\e{B}^+f(\r{a}[n];\nu_V),\leqno(5.32)
$$ 
provided $(3.8)$. }
\medskip 
We now turn to the contribution of Eisenstein series. We
see readily that
$$ 
e^{(1)}_p(\e{P}_0f;\nu)= {\zeta(z_2+\nu)\zeta(z_2-\nu)
\over\zeta(1-2\nu)}\left(\Phi_p^++\Phi_p^-\right)f(\nu).\leqno(5.33)
$$ 
This and $(5.6)$ confirm our claim on the convergence of
$(3.39)$ that is made prior to $(3.40)$. The discussion of
$\e{E}^{(1)}(\e{P}_0f;\r{a}[y],\nu)$ is obviously analogous to that
of
$\varpi_V\e{P}_0f(\r{a}[y])$ with $V$ in the unitary principal
series. Hence, it suffices to state
\medskip
\noindent {\bf Lemma 8.} {\it Let
$Z(\nu)=\zeta(z_1+\nu)\zeta(z_1-\nu)\zeta(z_2+\nu)\zeta(z_2-\nu)$.
Then
$$
\e{T}\varpi_E^{(1)}\e{P}_0f(1)=\int_{(0)} {Z(\nu)
\over\zeta(1+2\nu)\zeta(1-2\nu)}\left\{
\sum_{n=1}^\infty{\sigma_{2\nu}(n)\over n^{{1\over2}+\nu}}
\left(\e{B}^++\e{B}^-\right)f(\r{a}[n];\nu)\right\}{d\nu\over2\pi i},
\leqno(5.34)
$$ 
provided $(3.8)$. }
\medskip 
As to $\e{E}^{(0)}(\e{P}_0f;\r{a}[y],\nu)$, we observe that
the functional equation $(3.33)$ implies the relation
$d_p(\nu)e_p(\cdot;\nu)=e_p(\cdot;-\nu)$. Thus
$$
\e{E}^{(0)}(\e{P}_0f;1,\nu)=
\e{D}(\e{P}_0f;\nu)+\e{D}(\e{P}_0f;-\nu),\leqno(5.35)
$$ 
where
$$
\e{D}(\e{P}_0f;\nu)= {\zeta(z_2+\nu)\zeta(z_2-\nu)
\over\zeta(1-2\nu)}\left(\e{C}^++\e{C}^-\right)f(\nu),\leqno(5.36)
$$ 
with
$$
\e{C}^\delta f(\nu)=\sum_{p=-\infty}^\infty
\Phi_p^\delta f(\nu).\leqno(5.37)
$$
\medskip 
To compute this, we put
$$
\e{C}_\delta f(\nu)=\sum_{p=-\infty}^\infty\int_G  
f(\r{g})\e{A}^{-\delta}
\phi_{-p}(\r{g};-\nu)d\r{g},\leqno(5.38)
$$ 
which is regular for $|\Re\nu|<{1\over2}$, since the integral
satisfies the same bound as $(5.6)$. We have
$\e{C}_\delta f(\nu)=\e{C}^\delta f(\nu)$ on the imaginary axis. Let
us suppose $-{1\over2}<\Re\nu<0$. In $(5.38)$, use the first line of
$(2.16)$, but with the contour $\Im\xi={1\over2}$, so that the
quadruple integral converges absolutely; here we need
$(3.16)$. Take the integral over
$K$ innermost, and apply integration in parts many times, while
noting that $\partial_\theta^q f$ with any fixed $q$ still satisfies
the bound $(3.16)$. We see now that the sum over $p$ can be taken
inside the first triple integral. Then we may shift the
$\xi$-contour back to $\B{R}$. Undoing integration in parts, we get
$$
\leqalignno{
\e{C}_\delta f(\nu)&=\int_0^\infty
y^{-{3\over2}+\nu}\int_{-\infty}^\infty
\int_{-\infty}^\infty{\exp(2\pi i(y\xi-\delta x))\over
(\xi^2+1)^{-\nu+{1\over2}}}&(5.39)\cr
&\times\sum_{p=-\infty}^\infty\left({\xi+i\over\xi-i}\right)^{\delta
p}\int_0^\pi
f(\r{n}[x]\r{a}[y]\r{k}[\theta])\exp(-2pi\theta){d\theta\over\pi}\,
d\xi\,dx\,dy,  }
$$ 
and thus
$$
\leqalignno{
\e{C}_\delta f(\nu)&=\int_0^\infty y^{z_2-1+\nu}\int_{-\infty}^\infty
\int_{-\infty}^\infty{\exp(2\pi iy(\xi-\delta x))\over
(\xi^2+1)^{-\nu+{1\over2}}}f(\r{n}[x]\r{k}_\xi)
d\xi\,dx\,dy,\quad&(5.40) }
$$ 
with
$$
\r{k}_\xi={1\over\sqrt{\xi^2+1}}
\left[\matrix{\phantom{-}\xi&\delta\cr -\delta&\xi}\right]
\in K.\leqno(5.41)
$$ The last double integral over $(\xi,x)$ converges absolutely. We
take the $x$-integral innermost, and perform the change of variable
$x\mapsto \delta\xi+\delta x\left(\xi+\xi^{-1}\right)$, getting
$$
\leqalignno{
\qquad\e{C}_\delta f(\nu)&=\int_0^\infty y^{z_2-1+\nu}
\int_{-\infty}^\infty|\xi|^{w_1-w_4-1}
(\xi^2+1)^{z_2+\nu}\hat{\psi}_\tau
\left(\left(\xi+{1\over\xi}\right)y\right)\eta(\xi)d\xi\,dy&(5.42)\cr
&=\left((\hat{\psi}_\tau)^++(\hat{\psi}_\tau)^-\right)
\left(z_2+\nu\right)
\eta^*\left(z_3+\nu\right), }
$$ 
where $(\hat{\psi}_\tau)^\pm$ and $\eta^*$ are Mellin transforms
on
$(0,\infty)$ of $\hat{\psi}_\tau(\pm\cdot)$ and $\eta$, respectively;
and
$$ 
z_3={1\over2}(w_1+w_3-w_2-w_4+1).\leqno(5.43)
$$ 
In deducing $(5.42)$, we have used $(3.6)$ and $(5.21)$. The
second line of $(5.42)$ is a regular function of $\nu$ in a
neighbourhood of the imaginary axis; thus, it is equal to
$\e{C}^\delta f(\nu)$ if $\nu\in i\B{R}$.
\medskip 
Summing up, we obtain
\medskip
\noindent 
{\bf Lemma 9.} {\it
$$
\e{T}\varpi_E^{(0)}\e{P}_0f(1) =\int_{(0)}{Z(\nu)\over\zeta(1-2\nu)}
\left((\hat{\psi}_\tau)^++(\hat{\psi}_\tau)^-\right)
\left(z_2+\nu\right)
\eta^*\left(z_3+\nu\right){d\nu\over\pi i},\leqno(5.44)
$$ 
provided $(3.8)$. }
\vskip 1cm
\centerline{\bf 6. Limit}
\bigskip 
Collecting $(3.30)$, $(3.40)$, $(3.41)$, and Lemmas 6--9,
we arrive at a spectral decomposition of $\e{T}\e{P}f(1)$. Thus a
major part of our   proof of Theorem B has been finished. It remains
for us to discuss the behaviour   of this decomposition as $\eta$
tends to the characteristic function of
$\B{R}^\times$, and $\tau$ to that of negative reals. This limiting
procedure is our main task in the present section. Our   basic
implement here is the Mellin transform.
\medskip 
To facilitate various convergence issues, we restrict
$w$ by
$$
\left\{w:\,\hbox{bounded and $\Re w_3>\Re
w_2+3>4,\,\displaystyle{3\over2}+\Re w_1>
\Re w_4>\Re w_1+1>2$}\right\},\leqno(6.1)
$$ 
which is obviously contained in $(3.8)$. Also we set
$$
\eta(x)=\exp\left(- 
{1\over2}\xi_1\left(|x|+{1\over|x|}\right)\right),\quad
x\in\B{R}^\times,
\leqno(6.2)
$$ 
and
$$
\tau(x)=0,\,x\ge0;\quad
\tau(x)=\exp\left(- 
{1\over2}\xi_2\left(|x|+{1\over|x|}\right)\right),\,x<0,
\leqno(6.3)
$$ 
where both $\xi_1,\,\xi_2>0$ are supposed to tend to $0$. It is
immediate that $(3.9)$ implies
$$ 
{1\over2}\lim_{\xi_2\to0^+}\,
\lim_{\xi_1\to0^+}\e{T}\e{P} f(1)=\sum_{n=1}^\infty\sum_{m=1}^\infty
{\sigma_{w_1-w_4}(m)\sigma_{w_2-w_3}(m+n)\over
m^{w_1}(m+n)^{w_2}}\psi\left({m\over m+n}\right).\leqno(6.4)
$$ 
Via this relation a spectral decomposition of the right side
member arises.
\medskip 
We consider first the contribution of the unitary principal
series representations.  To this end, we transform $(5.24)$ into a
Mellin inverse integral: Put
$$
\psi_\tau^*(s)=\int_0^\infty {1\over (x+1)^{w_2}}\psi\left({1\over
x+1}\right)\tau(-x)x^{s-1} dx,\leqno(6.5)
$$ which is regular and of fast decay in any vertical strip in the
$s$-plane.
\medskip
\noindent 
{\bf Lemma 10.} {\it We have, for $\nu\in i{\Bbb R}$,
$$
\e{B}^\delta f(\r{a}[y];\nu)=y^{-z_3+{1\over2}}\int_{(0)}\eta^*(s_1)
\Psi_\tau^\delta(s_1;\nu)(2\pi y)^{-s_1}ds_1,\leqno(6.6)
$$ 
where $\eta^*$ is as in $(5.42)$, and
$$
\leqalignno{\qquad
\Psi_\tau^\delta(s_1;\nu)=&-4(2\pi)^{w_2-w_3-3}\int_{(\alpha)}
\psi_\tau^*(s_2)\cos\left({1\over2}\pi\big((1+\delta) (z_1-s_2)
+(1-\delta)\nu\big)\right)&(6.7)\cr
&\times\cos\left({1\over2}\pi\big (s_1+s_2+1-w_2-w_4-\delta
(s_2+1-w_1-w_2)\big)\right)\cr &\times\Gamma\left(z_1-s_2+\nu\right)
\Gamma\left(z_1-s_2-\nu\right)\Gamma(s_2+1-w_1-w_2)\cr
&\times\Gamma(s_1+s_2+1-w_2-w_4) ds_2. }
$$ 
Here $\alpha$ is to satisfy $\Re z_1>\alpha>\Re(w_1+w_2)-1$. Such
an  
$\alpha$ exists if $(6.1)$ holds. }
\medskip
\noindent 
{\it Proof.\/} We first transform $(5.20)$ into a Mellin
inverse   integral. To this end, we invert $(6.5)$ and have
$$ 
{1\over (x+1)^{w_2}}\psi\left({x\over
x+1}\right)\tau(-x)={1\over2\pi i}\int_{(\alpha)}
\psi_\tau^*(s)x^{s-w_2}ds,\leqno(6.8)
$$ 
with an arbitrary $\alpha$. Here we have used $\tau(x)=\tau(1/x)$.
Multiply both sides by the factor
$x^{-w_1}\exp(-ax-2\pi iux)$, $a>0$, and integrate over
$(0,\infty)$. The left side converges uniformly for $a\ge0$. Moving
the contour $(\alpha)$ to the right if necessary, the double integral
on the right side is seen to converge absolutely, provided $a>0$.  
Exchange the order of integration, compute the inner integral
explicitly, and   observe that the resulting integral is uniformly
convergent for $a\ge0$, because of the fast decay of $\psi_\tau^*$.
Thus we have, for $u\in{\Bbb   R}^\times$,
$$
\leqalignno{
\hat{\psi}_\tau(u)=&{1\over2\pi
i}\int_{(\alpha)}\psi_\tau^*(s)(2\pi|u|)^{w_1+w_2-1-s}
\Gamma(s+1-w_1-w_2)&(6.9)\cr &\times\exp\left(-{1\over2}\pi
i\,\sgn(u)(s+1-w_1-w_2)\right) ds, }
$$ with any $\alpha>\Re (w_1+w_2)-1$.
\par To the inner integral of $(5.24)$ we apply a similar procedure:
Multiply the integrand by the factor $\exp(-a|x|)$, $a>0$, replace
$\hat{\psi}_\tau$ by $(6.9)$, and exchange the order of integration.
We get the expression
$$
\leqalignno{ &{1\over2\pi i}\int_{(\alpha)}\psi_\tau^*(s_2) (2\pi
u)^{w_1+w_2-1-s_2}
\Gamma(s_2+1-w_1-w_2)&(6.10)\cr
\times\sum_\pm
\exp&\left(\pm{1\over2}\pi i\delta(s_2+1-w_1-w_2)\right)
\int_0^\infty x^{s_2-w_2-w_4}\exp(-(a\pm2\pi iy)x)\eta(x)dx\,ds_2. }
$$ 
On noting that $\eta(x)=\eta(1/x)$, use the Mellin inverse of
$\eta^*$. Because of the uniform convergence for
$a\ge0$, we see that the integral in question is equal to
$$
\leqalignno{ -{1\over2\pi^2}&\int_{(0)}\eta^*(s_1)
\int_{(\alpha)}\psi_\tau^*(s_2) (2\pi u)^{w_1+w_2-1-s_2}(2\pi
y)^{w_2+w_4-1-s_1-s_2}&(6.11)\cr
&\times\Gamma(s_2+1-w_1-w_2)\Gamma(s_1+s_2+1-w_2-w_4)\cr
&\times\cos\left({1\over2}\pi (s_1+s_2+1-w_2-w_4-\delta
(s_2+1-w_1-w_2))\right)ds_2ds_1, }
$$ 
provided $(6.1)$. Let us assume temporarily that $\alpha$ is such
that $0<\Re z_1-\alpha <{1\over4}$ as well as $\alpha>\Re
(w_1+w_2)-1$,   which does not conflict with $(6.1)$. We insert
$(6.11)$ into $(5.24)$. The resulting triple integral converges
absolutely, because of $(2.45)$. We   take the $u$-integral
innermost and invoke $(2.39)$--$(2.40)$ on noting we have presently
$\nu\in i{\Bbb R}$. In the result we can drop the condition $\Re
z_1-\alpha<{1\over4}$. This ends the proof.
\medskip The representation $(6.7)$ shows that
$\Psi_\tau^\delta(s_1;\nu)$ is in fact regular with respect to $\nu$
in a neighbourhood of the imaginary   axis. A shift to the far right
of the contour in
$(6.7)$ gives
$$
\Psi_\tau^\delta(s_1;\nu)\ll ((|s_1|+1)/(|\nu|+1))^C,\leqno(6.12)
$$ 
uniformly as $|\Re\nu|,\,|\Re s_1|<\varepsilon$, provided $(6.1)$
holds.  This is a consequence of the fast decay of $\psi_\tau^*$.
\medskip Let us now take $(6.2)$--$(6.3)$ into account precisely,
and consider   the limit of $\e{B}^\delta(\r{a}[y];\nu)$ as
$\xi_1,\,\xi_2$ tend to $0$: We claim first that
$$
\e{B}^\delta f(\r{a}[y];\nu)=2\pi iy^{-z_3+{1\over2}}
\Psi_\tau^\delta(0;\nu) +O\left(y^{-\Re z_3+{1\over2}}(\xi_1
y)^\varepsilon (|\nu|+1)^{-C}\right),\leqno(6.13)
$$ 
uniformly as $\xi_1\to0^+$. To confirm this, we shift the contour
in $(6.6)$ to $(\varepsilon)$, and note that, since $(6.2)$ implies
$\eta^*(s_1)=2K_{s_1}(\xi_1)$,
$$
\leqalignno{
\eta^*(s_1)&={\pi\over\sin\pi
s_1}\left(I_{-s_1}(\xi_1)-I_{s_1}(\xi_1)\right)&(6.14)\cr
&={\pi\over\sin\pi s_1} {1\over\Gamma(1-s_1)}
\left({\xi_1\over2}\right)^{-s_1} +O\left(\xi_1^{\Re
s_1}\exp(-|s_1|)\right), }
$$ 
with the implied constant being absolute. The bound
$(6.12)$ yields that the error-term contributes negligibly; and as to
the main term it suffices to shift the contour to
$(-\varepsilon)$.
\par 
We insert $(6.13)$ into $(5.25)$, and get
$$
\leqalignno{
\e{T}\varpi_V\e{P}_0f(1) &=2\pi
i\left(|\varrho_V(1)|^2+|\varrho_V(-1)|^2\right) H_V(z_1)
H_V(z_2)H_V(z_3)&(6.15)\cr &\times\left(\Psi_\tau^+
+\epsilon_V\Psi_\tau^-\right)(0;\nu_V)
+O(|\varrho_V(1)|^2\xi_1^\varepsilon(|\nu_V|+1)^{-C}), }
$$ 
in which we have used the fact that $(6.1)$ implies $\Re
z_3>{5\over4}$. Because of $(2.29)$, we find that
$$
\leqalignno{
\lim_{\xi_1\to0^+}\sum_V
\e{T}\varpi_V\e{P}_0f(1)& =2\pi
i\sum_V\left(|\varrho_V(1)|^2+|\varrho_V(-1)|^2\right)&(6.16)\cr
&\times H_V(z_1)H_V(z_2)H_V(z_3)\left(\Psi_\tau^+
+\epsilon_V\Psi_\tau^-\right)(0;\nu_V), }
$$ 
with $V$ running over all irreducible representations in the
unitary principal series.
\par 
Next, we observe that for $\Re s>0$
$$
\psi^*_\tau(s)={1\over\pi i}\int_{(0)}K_\mu(\xi_2)
\psi^\ast(s-\mu)d\mu.\leqno(6.17)
$$ 
Here $\psi^*$ is defined to be the right side of $(6.5)$ without
the factor $\tau(-x)$. It is regular and of fast decay for $\Re
s>0$. We have
$$
\psi^*_\tau(s)=I_0(\xi_2)\psi^*(s)+\psi^{**}_\tau(s),\leqno(6.18)
$$ 
with
$$
\psi^{**}_\tau(s)=-{1\over2i}\int_{(\varepsilon)}
I_\mu(\xi_2)\left(\psi^*(s+\mu)+
\psi^*(s-\mu)\right){d\mu\over\sin\pi\mu}.\leqno(6.19)
$$ 
For $\Re s>0$ the $\psi^{**}_\tau(s)$ is regular and
$\ll \xi_2^\varepsilon(|s|+1)^{-C}$. This gives readily
$$
\Psi_\tau^\delta(0;\nu)=\Psi^\delta(\nu)+O(\xi_2^\varepsilon
(|\nu|+1)^{-C}),\leqno(6.20)
$$ 
where $\Psi^\delta(\nu)$ is defined by $(6.7)$ with $s_1=0$ and
$\psi^*$ in place of $\psi^*_\tau$.
\medskip 
Hence, we have proved, as a consequence of Lemmas 6 and 10,
\medskip
\noindent 
{\bf Lemma 11.} {\it Let $(6.1)$ hold. Then, with $V$
running over all irreducible representations in the unitary
principal series, we have
$$
\leqalignno{ {1\over2}\lim_{\xi_2\to0^+}\lim_{\xi_1\to
0^+}\sum_V\e{T}&\varpi_V\e{P}_0f(1) =\pi
i\sum_V\left(|\varrho_V(1)|^2+|\varrho_V(-1)|^2\right) &(6.21) \cr
&\times H_V(z_1)H_V(z_2)H_V(z_3)\left(\Psi^+
+\epsilon_V\Psi^-\right)(\nu_V). }
$$ 
Here
$$
\leqalignno{\qquad
\Psi^\delta(\nu)=&-4(2\pi)^{w_2-w_3-3}\int_{(\alpha)}
\psi^*(s)\cos\left({1\over2}\pi\big((1+\delta)(z_1-s)
+(1-\delta)\nu\big)\right)&(6.22)\cr
&\times\cos\left({1\over2}\pi\big (s+1-w_2-w_4-\delta
(s+1-w_1-w_2)\big)\right)\cr &\times\Gamma\left(z_1-s+\nu\right)
\Gamma\left(z_1-s-\nu\right)\Gamma(s+1-w_1-w_2)\cr
&\times\Gamma(s+1-w_2-w_4)ds, }
$$ 
with $\Re z_1>\alpha>\Re(w_1+w_2)-1$ and
$$
\psi^*(s)=\int_0^\infty {x^{s-1}\over (x+1)^{w_2}}\psi\left({1\over
x+1}\right)dx.\leqno(6.23)
$$
}
\medskip
We compare $(6.21)$ with
$(4.5.5)$ and $(4.5.9)$ of [19], and $(6.22)$ with  
$(4.4.16)$--$(4.4.17)$ there. To facilitate it, we give the table:
$$
\eqalign{ (w_1,w_2,w_3,w_4)&\mapsto (u,w,z,v),\cr
|\varrho_V(\pm1)|^2&\mapsto {1\over4}\alpha_j,\cr
\Psi^\delta&\mapsto {2\over\pi i}\Phi_\delta,\cr
\psi^*&\mapsto \tilde{g}, }\leqno(6.24)
$$ 
where on the left are our present objects and on the right those
corresponding in Sections 4.3--4.4 of [19]. Likewise 
$\epsilon_V\mapsto
\epsilon_j$, $H_V\mapsto H_j$. Note that
$(4.4.16)$ in [19] is to be corrected: the second $\xi$ on the right
side should have the opposite sign. We find that the agreement is
perfect. This ends the treatment of the unitary principal series.
\medskip Concerning the contribution of the discrete series, we
return to Lemma 7. We observe that Lemma 10 holds for
$\e{B}^+\left(\r{a}[y];\ell-{1\over2}\right)$ with ${\Bbb  
Z}\ni\ell\ge1$ as well. In fact, the necessary change in the proof
takes place only after
$(6.11)$. The use of $(2.39)$--$(2.40)$ is replaced by that of  
$(2.43)$. This and the argument leading to Lemma 11 yield the
following assertion,   which in view of $(2.20)$ and $(2.28)$ 
coincides with
$(4.5.6)$ of [19]:
\medskip
\noindent 
{\bf Lemma 12.} {\it Let $(6.1)$ hold. Then, with $V$
running over all irreducible representations in the discrete series,
we have
$$
\leqalignno{ {1\over2}&\lim_{\xi_2\to 0^+}\lim_{\xi_1\to
0^+}\sum_V\e{T}\varpi_V\e{P}_0f(1)&(6.25) \cr &=\pi
i\sum_V\left(|\varrho_V(1)|^2+|\varrho_V(-1)|^2\right) H_V(z_1)H_V
(z_2)H_V(z_3)\Psi^+(\nu_V), }
$$ 
where
$$
\leqalignno{ &\Psi^+(\nu_V)=2(-1)^{\ell-1}(2\pi)^{w_2-w_3-2}
\cos\left(\txt{1\over2}\pi(w_1-w_4)\right)&(6.26)\cr
&\times\int_{(\alpha)}
\psi^*(s)
{\Gamma(\ell-{1\over2}+z_1-s)\over\Gamma(\ell+{1\over2}-z_1+s)}
\Gamma(s+1-w_1-w_2)\Gamma(s+1-w_2-w_4)ds. }
$$ 
with $\nu_V=\ell-{1\over2}$ and $\alpha$ as above. }
\medskip
\noindent 
{\bf Remark.} With respect to the process behind the
appearance of products of three values of Hecke series, there arises
a   notable difference between the present work and [19]. Here
$H_V(z_1)$ comes from
$(2.26)$, $H_V(z_2)$ from $(3.7)$, and
$H_V(z_3)$ from the above limiting procedure. In [19] the
$H_V(z_1)$ corresponds to the sum over the shift parameter, i.e.,
the  $n$ of $(3.9)$ and thus in a context similar to the
present. However,   there
$H_V(z_2)$ and $H_V(z_3)$ appear combined in the product
$H(z_2)H(z_3)$   as a consequence of the multiplicativity of Hecke
operators that is   irrelevant to our argument (but $(2.30)$ depends
on it via the proof). This   observation applies also to the product
of six values of the zeta-function in the numerator of
$(6.27)$ below.
\medskip 
As to the consequence of Lemma 8, it should be enough to
state only the end result, which coincides with the sum of $(4.5.4)$
and $(4.5.8)$ of   [19]:
\medskip
\noindent 
{\bf Lemma 13.} {\it Let $(6.1)$ hold. Then we have
$$
\leqalignno{ &\qquad{1\over2}\lim_{\xi_2\to0^+}\lim_{\xi_1\to0^+}
\e{T}\varpi_E^{(1)}\e{P}_0f(1)&(6.27)\cr &={1\over2}\int_{(0)}
{\zeta(z_1+\nu)\zeta(z_1-\nu)\zeta(z_2+\nu)\zeta(z_2-\nu)
\zeta(z_3+\nu)\zeta(z_3-\nu)
\over\zeta(1+2\nu)\zeta(1-2\nu)}
\left(\Psi^++\Psi^-\right)(\nu){d\nu}, }
$$ 
}
\medskip 
It remains to discuss $\e{T}\e{P}\!{}_\infty f(1)$ and
$\e{T}\varpi_E^{(0)}\e{P}_0f(1)$. The following assertion   coincides
with $(4.3.16)$ of [19]:
\medskip
\noindent 
{\bf Lemma 14.} {\it Let $(6.1)$ hold. Then we have
$$
\leqalignno{ {1\over2}&\lim_{\xi_2\to0^+}
\lim_{\xi_1\to0^+}\e{T}\e{P}\!{}_\infty
f(1)=\psi^*(w_1+w_2-1)&(6.28)\cr
&\times{1\over\zeta(2z_2+1)}\zeta(w_1+w_2-1)\zeta(w_3+w_4)
\zeta(w_3-w_2+1)\zeta(w_4-w_1+1), }
$$ 
and
$$
\leqalignno{ {1\over2}&\lim_{\xi_2\to0^+}
\lim_{\xi_1\to0^+}\e{T}\varpi_E^{(0)}\e{P}_0f(1)
=\psi^*(w_2+w_4-1)&(6.29)\cr
&\times{1\over\zeta(2z_3+1)}\zeta(w_1+w_3)\zeta(w_2+w_4-1)
\zeta(w_3-w_2+1)\zeta(w_1-w_4+1). }
$$ 
}
\medskip
\noindent 
{\it Proof.\/} To confirm the former, it suffices to
observe $(3.30)$, the first equation in $(3.41)$, and  that $(6.5)$
gives
$\hat{\psi}_\tau(0)=\psi_\tau^*(w_1+w_2-1)$. As to the latter, we
apply, on the right side of $(5.44)$, the change of variable
$\nu\mapsto \nu-z_3$, and shift the contour to
$(\varepsilon)$. We do not encounter
any singularity, because of
$(6.1)$. Following the argument for $(6.13)$, we have
$$
\lim_{\xi_1\to0^+}\e{T}\varpi_E^{(0)}\e{P}_0f(1)
=2{Z(-z_3)\over\zeta(2z_3+1)}\left((\hat{\psi}_\tau)^+
+(\hat{\psi}_\tau)^-\right)(w_4-w_1).
\leqno(6.30)
$$ 
The inversion of the formula $(6.9)$ gives
$$
\leqalignno{
\left((\hat{\psi}_\tau)^+ +(\hat{\psi}_\tau)^-\right)&(w_4-w_1)
&(6.31)\cr =2(2\pi)^{w_1-w_4}&\cos\left(
{1\over2}\pi(w_1-w_4)\right)\Gamma(w_4-w_1)
\psi_\tau^*(w_2+w_4-1). }
$$ 
This and the functional equation for the zeta-function yield
$(6.29)$.
\medskip
\noindent 
{\bf Remark.} In [19] those two terms corresponding to
$(6.28)$ and   $(6.29)$ come up as residual terms. Here they are,
respectively, understood as   the consequence of the automorphic
regularization $(3.24)$ and the   contribution of the constant terms
of Eisenstein series.
\bigskip 
Gathering $(3.40)$, $(6.4)$ and Lemmas 11--14 together, we
see that we   have now proved $(4.3.15)$--$(4.3.16)$ and Lemmas 4.5
and 4.6 of [19], with a new argument. One should note that the
domain $(6.1)$ is   not disjoint with $(4.3.10)$ of [19]. Those
assertions of [19] give rise to   the spectral decomposition of
$J(w;g)$ defined at $(1.8)$. Hence, we have achieved the same. What
remains is to continue analytically the spectral decomposition of  
$J(w;g)$ to a neighbourhood of $\r{p}_{1\over2}$. This is naturally
the same as Sections 4.6--4.7 of [19], and we skip it.
\medskip Thus, for instance, the sum $(6.21)$ converges absolutely
uniformly in a neighbourhood of $\r{p}_{1\over2}$. The value at
$\r{p}_{1\over2}$ is equal to
$$ 
{1\over2}\sum_V\left(|\varrho_V(1)|^2+|\varrho_V(-1)|^2\right)
H_V\left({1\over2}\right)^3
\left(\psi_++\psi_-\right)(\nu_V),\leqno(6.32)
$$ 
where $V$ are in the unitary principal series, and
$\psi_\delta(\nu)$   is the value of $2\pi i\Psi^\delta(\nu)$ with
$w=\r{p}_{1\over2}$. Here we   have used the fact that
$H_V\left({1\over2}\right)=0$ if $\epsilon_V=-1$, as   implied by
$(2.32)$. The formula
$(6.22)$ gives, for $\nu\in i{\Bbb R}$,
$$
\leqalignno{
\psi_\delta(\nu)&={1\over\pi^2i}\int_{({1\over4})}
\psi^*(s)\cos\left({1\over2}\pi
\left((1+\delta)\left({1\over2}-s\right)+(1-\delta)\nu\right)\right)
&(6.33)\cr &\times\cos\left({1\over2}\pi(1-\delta)s\right)
\Gamma\left({1\over2}-s+\nu\right)
\Gamma\left({1\over2}-s-\nu\right)\Gamma(s)^2ds. }
$$ 
Inserting $(6.23)$ with $w_2={1\over2}$, we get an absolutely
convergent double integral. Exchanging the order of integration and
invoking
$(2.41)$--$(2.42)$, we find that
$$
\psi_\delta(\nu)=\int_0^\infty {1\over x(x+1)^{1\over2}}
\psi\left({1\over x+1}\right)\int_0^\infty j_0(-\delta
u)j_\nu(\delta   xu) {du\over u^{3\over2}}dx,\leqno(6.34)
$$ 
which implies $(2.34)$ with an obvious specialization of $\psi$.
The real reason for the appearance of $j_\nu$ here should of course
be   traced back to $(4.9)$. This concerns the contribution of
unitary principal   series representations only, but other parts are
dealt with similarly.   Therefore we have fully proved Theorem B,
with a method based solely upon the   spectral theory of
$L^2(\varGamma\backslash G)$.
\bigskip
\noindent 
{\bf Concluding Remark.} Theorem A is not an isolated
fact in the theory of zeta-functions. It has been extended to
$\e{M}_4(\zeta_F;g)$, where $\zeta_F$ are the Dedekind
zeta-functions of certain quadratic number fields $F$ with class
number one. See [5] and [7] (also  [6, Part X]). The relevant Lie
groups are $\r{PSL}_2(\B{C})$ and the product of two copies of
$\r{PSL}_2(\B{R})$, respectively, for the imaginary and the real
quadratic cases. What is interesting is that these extensions of
$\e{M}_4(\zeta;g)$ admit formulations highly analogous to Theorem B,
despite the difference of the representation theories of the two Lie
groups from that of $\r{PSL}_2(\B{R})$.  The construction
$(2.37)$--$(2.38)$ extends with cubic powers of central values of  
analogous Hecke series. Besides, $(2.34)$--$(2.35)$ extends with
Bessel functions of representations of corresponding Lie groups.
\par 
Also, an analogue of Theorem A is known for $\e{M}_2(H;g)$ with
Hecke series $H$ attached to any   holomorphic cusp form. Again a
formulation similar to Theorem B is possible. See   [18] and [20].
Hence, Theorem B could be a typical instance of a certain   general
structure among mean values of automorphic $L$-functions. However,
the   case of real analytic cusp forms is yet to be included in this
general   picture. Moreover, it is not known if the argument of the
present work, i.e., the Poincar\'e series approach, extends to
$\e{M}_2(H;g)$. To these topics  we shall return elsewhere.
\par Finally, it should be stressed that M. Jutila has developed a
method with which one can treat $L$-functions of all types equally
as far as   the underlying group is $\r{PSL}_2({\Bbb R})$, although
it does not produce results as explicit as Theorem A. See his
important work [14].
\bigskip
\noindent 
{\it Acknowledgement.\/} This work is an outcome of our
stay, in March and June 2002, at the Max Planck Institute for
Mathematics at Bonn, in the Special Activity in Analytic Number
Theory. We thank the MPIM and the organizers of the activity for the
invitation and the hospitality they have shown us.
\vskip 1cm
\centerline{\bf References}
\bigskip
\item{[1]}  {\it F.V. Atkinson\/}, The mean value of the Riemann
zeta-function, Acta Math.\ {\bf 81} (1949), 353--376.
\item{[2]} {\it E.M. Baruch\/} and {\it Z. Mao\/}, Bessel identities
in Waldspurger correspondence, the archi\-medean theory,  to appear
in Israel J. Math.
\item{[3]} {\it R.W. Bruggeman\/}, Fourier coefficients of
automorphic forms, Lecture Notes in Math., {\bf 865},
Springer-Verlag, Berlin 1981.
\item{[4]} {\it R.W. Bruggeman\/}, Automorphic forms,  Banach Center
Publ.\ {\bf 17} (1985), 31--74.
\item{[5]} {\it R.W. Bruggeman\/} and {\it Y. Motohashi\/}, Fourth
power moment of Dedekind zeta-functions of real quadratic number
fields with   class number one,  Functiones et Approximatio {\bf 29}
(2001), 41--79.
\item {[6]} {\it R.W. Bruggeman\/} and {\it Y. Motohashi\/}, A note
on   the mean value of the zeta and $L$-functions, X, Proc.\ Japan
Acad.\ {\bf 77A} (2001), 111--114;\ XIII, ibid {\bf 78A} (2002),
87--91.
\item{[7]} {\it R.W. Bruggeman\/} and {\it Y. Motohashi\/}, Sum
formula   for Kloosterman sums and the fourth moment of the Dedekind
zeta-function   over the Gaussian number field, Functiones et
Approximatio {\bf 31} (2003),   7--76.
\item{[8]} {\it J.W. Cogdell\/} and {\it I. Piatetski-Shapiro\/}, The
arithmetic and spectral analysis of Poin\-car\'e series,
Perspectives in Math.\ {\bf 13}, Academic Press, San Diego 1990.
\item{[9]} {\it J. Dixmier\/} and {\it P. Malliavin\/}, 
Factorisations   de fonction et de vecteurs ind\'efiniment
diff\'erentiables, Bull.\ Soc.\   Math.\ France {\bf 102} (1978),
305--330.
\item{[10]} {\it I.M. Gel'fand\/}, {\it M.I. Graev\/} and {\it I.I.
Pyatetski-Shapiro\/},  Representation theory and automorphic
functions, W.B. Saunders Company, Philadelphia 1969.
\item{[11]} {\it I.S. Gradshteyn\/} and {\it I.M. Ryzhik\/}, Tables
of integrals, series, and products,  Academic Press, San Diego 1979.
\item{[12]} {\it A. Ivi\'c\/}, Mean values of the Riemann  
zeta-function, Tata IFR Lect.\ Math.\ Phys.\ {\bf 82},
Springer-Verlag, Berlin-Heidelberg-New York-Tokyo 1991.
\item{[13]} {\it A. Ivi\'c\/}, On the error term for the fourth
moment   of the Riemann zeta-function,  J. London Math.\ Soc.\ {\bf
60} (1999), 21--32.
\item {[14]} {\it M. Jutila\/}, Mean values of Dirichlet series via  
Laplace transforms, In: Analytic Number Theory, Proc.\ 39th
Taniguchi Intern.\ Symp.\ Math., Kyoto 1996, ed.\ Y. Motohashi,
Cambridge Univ.\ Press, Cambridge 1997, pp.\ 169--207.
\item{[15]} {\it A.A. Kirillov\/},  On $\infty$-dimensional unitary
representations of the group of second-order matrices with elements
from a locally compact field, Soviet Math.\ Dokl., {\bf 4} (1963),
748--752.
\item{[16]} {\it Y. Motohashi\/}, The fourth power mean of the
Riemann zeta-function. In: Proc.\ Conf.\ Analytic Number Theory,
Amalfi 1989, eds.\ E. Bombieri et al., Univ.\ di Salerno, Salerno
1992, pp.\ 325--344.
\item{[17]} {\it Y. Motohashi\/}, An explicit formula for the
fourth   power mean of the Riemann zeta-function, Acta Math. {\bf
170} (1993),   181--220.
\item{[18]}  {\it Y. Motohashi\/}, The mean square of Hecke
$L$-series attached to holomorphic cusp forms, RIMS Kyoto Univ.
Kokyuroku {\bf 886} (1994), 214--227.
\item{[19]} {\it Y. Motohashi\/}, Spectral theory of the Riemann
zeta-function, Cambridge Tracts in Math.\ {\bf 127}, Cambridge
Univ.\ Press, Cambridge 1997.
\item{[20]} {\it Y. Motohashi\/}, A note on the mean value of the
zeta   and
$L$-functions, XII, Proc.\ Japan Acad.\ {\bf 78A} (2002), 36--41.
\item{[21]} {\it N.Ja.\ Vilenkin\/}, Special functions and the
theory of group representations, Amer.\ Math.\ Soc., Providence 1968.
\item{[22]} {\it N.Ja.\ Vilenkin\/} and {\it A.U. Klimyk\/},  
Representations of Lie groups and special functions, Vol.\ 1, Kluwer
Acad.\ Publ., Dordrecht-Boston-London 1991.
\item{[23]} {\it G.N. Watson\/}, A treatise on the theory of Bessel
functions, Cambridge Univ.\ Press, Cambridge 1996.
\item{[24]} {\it E.T. Whittaker\/} and {\it G.N. Watson\/},  A
course of modern analysis, Cambridge Univ.\ Press, London 1969.

\bigskip
\noindent
\font\small=cmr8 {\small Roelof W. Bruggeman
\par\noindent Department of Mathematics, Utrecht University,
\par\noindent P.O.Box 80.010, TA 3508 Utrecht, the Netherlands
\par\noindent Email: bruggeman@math.uu.nl}
\medskip
\noindent {\small Yoichi Motohashi
\par\noindent Department of Mathematics, College of Science and
Technology, Nihon University,
\par\noindent Surugadai, Tokyo 101-8308, Japan
\par\noindent Email: ymoto@math.cst.nihon-u.ac.jp,
am8y-mths@asahi-net.or.jp }

\bye